\documentclass{amsart}
\usepackage{amssymb}
\usepackage{amscd}
\usepackage{amsmath}

\newcommand{\Q}{\mathbb Q}
\newcommand{\Z}{\mathbb Z}
\newcommand{\R}{\mathbb R}
\newcommand{\C}{\mathbb C}

\newcommand{\F}{\mathbb F}

\newcommand{\A}{\mathbb A}
\newcommand{\N}{\mathbb N}

\makeatletter
\def\l@section{\@tocline{1}{4pt}{1pc}{}{}}
\def\l@subsection{\@tocline{2}{0pt}{2pc}{5pc}{}}
\makeatother
\title{Modularity of solvable Artin representations of GO$(4)$-type}

\author{Dinakar Ramakrishnan}
\thanks{Partially supported by the NSF grant DMS-9801328}
\address{253-37 Caltech, Pasadena, CA 91125}
\email{dinakar@its.caltech.edu}

\begin{document}
\maketitle

\pagestyle{myheadings}
\markboth{Dinakar Ramakrishnan}{Solvable Artin representations of 
GO$(4)$-type}
\bigskip

\qquad \qquad \qquad \qquad \qquad \qquad \qquad \qquad \qquad{\it To my teacher, Herv\'e Jacquet}

\bigskip

\section*
{\bf Introduction}

\bigskip

Let $F$ be a number field, and $(\rho, V)$ a continuous, $n$-dimensional representation 
of the 
absolute Galois group Gal$(\overline F/F)$ on a finite-dimensional $\C$-vector space $V$. 
Denote by $L(s, \rho)$ the associated $L$-function, which is known to be meromorphic with 
a functional equation. Artin's conjecture predicts that 
$L(s, \rho)$ is holomorphic everywhere except possibly at $s=1$, where its order of pole 
is the multiplicity of the trivial representation in $V$. The {\it modularity conjecture} 
of Langlands for such representations ([La3]), 
often called the {\it strong Artin conjecture}, 
asserts that there should be an 
associated (isobaric) automorphic form $\pi = \pi_\infty \otimes \pi_f$ on GL$(n)/F$ 
such that 
$L(s, \rho) = L(s, \pi_f)$. Since $L(s, \pi_f)$ possesses the requisite properties 
([JS]), the
modularity conjecture implies the Artin conjecture.

\medskip

For any field $k$, let ${\text GO}(n, k)$ denote the subgroup of ${\text GL}(n, k)$ 
consisting of 
{\it orthogonal similitudes}, i.e., matrices $g$ such that ${}^tg g = \lambda_g I$, with 
$\lambda_g \in k^\ast$. 

\medskip

We will say that a $k$-representation $(\rho, V)$  is {\it of GO$(n)$-type} iff dim$(V) = n$ and 
it factors as
$$
\rho \, = \, [{\text Gal}(\overline F/F) \, \stackrel{\sigma}{\longrightarrow} \, {\text GO}(n, k) \, 
\subset \, {\text GL}(V)].
$$

\medskip

In this article we prove 

\medskip

\noindent{\bf Theorem A} \, \it Let $F$ be a number field and 
let $(\rho, V)$ be a continuous, $4$-dimensional $\C$-representation of 
Gal$(\overline F/F)$ whose image is solvable and lies in GO$(4, \C)$. Then 
$\rho$ is modular, i.e., $L(s, \rho) = L(s, \pi_f)$ for some isobaric automorphic representation 
$\pi = \pi_\infty \otimes \pi_f$ of GL$(4, \A_F)$. Moreover, $\pi$ is cuspidal iff
$\rho$ is irreducible.
\rm

\medskip

Among the finite groups $G$ showing up as the images of such $\rho$ are those fitting into an exact sequence
$$
1 \, \rightarrow \, C \, \rightarrow \, H \times H \, \rightarrow \, G \, \rightarrow \, \{\pm 1\} \, \rightarrow \, 1,
$$
where $H$ is any finite solvable group in GL$(2, \C)$ with scalar subgroup $C$, the embeding of $C$ in $H \times H$ is given by $x \to (x, x^{-1})$, and the action of $\{\pm 1\}$ on $(H \times H)/C$ is induced by the permutation of the two factors of $H \times H$. This is due to the well known fact (see section 1) that GO$(4, \C)$ contains a subgroup of index $2$ which is a quotient of GL$(2, \C) \times$ GL$(2, \C)$ by $\C^\ast$. Of particular interest are the examples where $H$ is a central extension of $S_4$ or $A_4$ (cf. section 8).  

One can ask if this helps furnish new examples of Artin's conjecture, and the answer is yes.

\medskip

\noindent{\bf Corollary B} \, \it Let $F$ be a number field, 
and let $\rho, \rho'$ be 
continuous $\C$-representations of Gal$(\overline F/F)$ of 
solvable GO$(4)$-type. Then
Artin's conjecture 
holds for $\rho \otimes \rho'$. 
\rm

\medskip
 
We will show in section 8 that in fact there is, for each $F$, a doubly {\it infinite} family of 
such examples where the representations $\rho \otimes \rho'$ are {\it irreducible} and 
{\it primitive} (i.e., not induced) of dimension $16$. 
Primitivity is important because 
Artin $L$-functions are inductive, and one wants to make sure that these examples do {\sl not} 
come by induction from known (solvable) cases in low dimensions. 
We will then show (in secton 9) that, given 
any $\rho$ as in Theorem A with corresponding extension $K/F$, the {\it strong Dedekind conjecture} 
holds for certain non-normal extensions $N/F$ contained in $K/F$, for instance when $[N:F] = 3^a$, 
$a \geq 1$ (see Theorem 9.3); the assertion is that the ratio $\zeta_N(s)/\zeta_F(s)$ is the 
standard $L$-function of an isobaric automorphic form $\eta$ on GL$(3^a-1)/F$. 
It implies that for any cusp form $\pi$ on GL$(n)/F$, the formally defined Euler product 
$L(s, \pi_N)$ (see the discussion before {\it Corollary 9.4}) admits a meromorphic continuation and functional equation, and more importantly, it is divisible by $L(s, \pi)$. It should be noted that we do not know if $\pi_N$ is automorphic.
For a curious consequence of this result see Remark 9.6. 

\medskip

It has been known for a long time, thanks to the results of Artin and Hecke, 
that monomial representations of Gal$(\overline F/F)$, 
i.e., those induced by one-dimensional representations of Gal$(\overline F/K)$ with $K/F$ finite, 
satisfy Artin's conjecture. (In fact this holds for any multiple of a monomial representation.) But the strong Artin conjecture is still open for these except when $K/F$ is 
normal and solvable ([AC]) and when $[K:\Q] = 3$ ([JPSS1]. The work [AC] of Arthur and Clozel implies that the 
strong conjecture holds for $\rho$ with nilpotent image, in fact whenever $\rho$ is {\it accessible} ([C]), 
i.e., a positive integral linear combination of representations induced from linear characters of open, 
subnormal subgroups. It should however be noted that Dade has shown ([Da]) that all accessible 
representations of solvable groups are monomial. 

\medskip
 
The odd dimensional orthogonal representations with solvable image are simpler than the even dimensional
ones. Indeed we have

\medskip

\noindent{\bf Proposition C} \, \it Let $\rho$ be a continuous, irreducible, solvable $\C$-representation 
of Gal$(\overline F/F)$ of GO$(n)$-type with $n$ odd. Then $\rho$ is monomial 
and hence satisfies the Artin conjecture. If $\rho$ is in addition self-dual, it must be induced 
by a quadratic character. 
\rm

\medskip

By the groundbreaking 
work of Langlands in the seventies ([La1]), and the ensuing theorem of Tunnell in 1980 ([Tu]), one knows 
that the strong 
Artin conjecture holds for all two-dimensional representations $\sigma$ with solvable image. 
The {\it raison d'etre} for this article is the desire to find irreducible, {\it solvable}, but 
non-monomial, even primitive, examples in higher dimensions satisfying the Artin cojecture. We were led to this problem a few years ago by a remark of J.-P.Serre. One could look at the 
symmetric powers of a two-dimensional {\it solvable} $\sigma$, but they are all 
in the linear span of 
monomial representations and $1$-dimensional twists of $\sigma$. The facts that the symmetric square 
of $\sigma$ is modular (by Gelbart-Jacquet [GeJ]) and monomial (for $\sigma$ solvable) were used 
crucially in [La1] and [Tu].

In the {\it non-solvable} direction, which is orthogonal (no pun intended) to the 
one pursued in this paper, 
there has been some spectacular progress recently. For {\it odd} $2$-dimensionals $\rho$ of Gal$(\overline \Q/\Q)$ with projectivization $\overline \rho$ 
of $A_5$-type, a theorem of Buzzard, Dickinson, Shepherd-Barron and 
Taylor ([BDST]) 
establishes the modularity conjecture assuming the following: (i) $\overline \rho$ is unramified at $2$ and $5$, and (ii) $\overline \rho({\rm Frob}_2)$ has order $3$. In a sequel ([T]), this is shown with different ramification conditions, namely that under $\overline \rho$ (i') the inertia group at $3$ has odd order, and (ii') the decomposition group at $5$ is unramified at $2$. See also 
[BS] for some 
explicit examples. Some positive examples were given earlier in [Bu], and then in [Fr]. 
Moreover, the very recent theorems of Kim and 
Shahidi ([KSh]), and 
Kim ([K]), establishing the automorphy of the symmetric cube, and the symmetric 
fourth power, on 
GL$(2)$, establishes the strong Artin conjecture for sym$^3(\sigma)$, and sym$^4(\sigma)$, 
for all the 
$\sigma$ proved modular by [BDST]. 

If $\sigma, \sigma'$ are Galois representations which are modular, then one can deduce the Artin 
conjecture for $\sigma \otimes \sigma'$ by applying the Rankin-Selberg theory on GL$(m) \times $GL$(n)$ 
developed in the works of Jacquet, Piatetski-Shapiro, Shalika ([JPSS2], [JS]) and of Shahidi ([Sh1,2]); see 
also [MW]. This explains why the strong form of Artin's conjecture is really a bit stronger than the original conjecture that Artin made, at least given what 
one knows today. 

If $\sigma, \sigma'$ are $2$-dimensional representations which are modular, then the strong Artin 
conjecture for $\sigma \otimes \sigma'$ 
follows from the main theorem of [Ra1], hence the Artin conjecture holds, by the remark above, 
for $4$-fold tensor products of such representations. Now let $K/F$ be a quadratic extension with 
non-trivial automorphism 
$\theta$, and let $\sigma^{\theta}$ denote the $\theta$-twisted representation defined by 
$x \to \sigma(\tilde \theta x \tilde \theta^{-1})$, where $\tilde \theta$ is any lift to 
Gal$(\overline F/F)$. (The equivalence class of $\sigma^{\theta}$ is independent of the 
choice of $\tilde \theta$.) Given any irreducible, non-monomial $2$-dimensional representation 
$\sigma$ of Gal$(\overline F/K$) 
which is not isomorphic to any one-dimensional twist of $\sigma^{\theta}$ (see section 3), 
there is an irreducible 
$4$-dimensional representation $As(\sigma)$ of Gal$(\overline F/F)$ whose restriction to Gal$(\overline F/K)$ is 
isomorphic to $\sigma \otimes \sigma^{\theta}$. When $\sigma$ is of solvable type, one can combine the theorem of
Langlands-Tunnell with that of Asai ([HLR]) to deduce the Artin conjecture for $As(\sigma)$. 

\medskip

One of the main steps in our proof of Theorem A is that even modularity can be 
established for any Asai representation $As(\sigma)$ (in the solvable case). To begin, there exists, 
by Langlands-Tunnell, 
a cusp form $\pi$ on GL$(2)/K$ such that $L(s, \sigma) = L(s, \pi)$. It follows that 
$L(s, \sigma \otimes \sigma^{\theta})$ equals the Rankin-Selberg $L$-function 
$L(s, \pi \times (\pi \circ \theta))$. By [Ra1], there exists an automorphic form 
$\pi \boxtimes (\pi \circ \theta)$ on GL$(4)/K$ such that $L(s, \sigma \otimes \sigma^{\theta}) = 
L(s, \pi \boxtimes (\pi \circ \theta))$. Since $\pi \boxtimes (\pi \circ \theta)$ is $\theta$-invariant, 
one can now deduce the existence of a cusp form $\Pi$ on GL$(4)/F$ whose base change to $K$ is 
$\pi \boxtimes (\pi \circ \theta)$, but even when $\pi \boxtimes (\pi \circ \theta)$
is cuspidal, $\Pi$ is unique only up to twisting by the quadratic character 
$\delta$, say, of the idele class group of $F$ corresponding to 
$K/F$ by class field theory. But it is not at all clear 
that $\Pi$, or $\Pi \otimes \delta$, should correspond precisely to $As(\sigma)$, with an identity of the 
corresponding $L$-functions. (It is easy to see that 
the local factors agree at half the places.) Put another way, one can construct an irreducible admissible 
representation $As(\pi)$ of GL$(4, \A_F)$ which has the same local factors as does $As(\sigma)$. But the 
problem 
is that there is no simple reason why $As(\pi)$ should be automorphic, {\it even} 
when $\pi$ is dihedral. (Recall that $\pi$ is dihedral, or of CM type, iff 
it is associated to a 
representation $\sigma$ 
of the global Weil group $W_K$ induced by a character $\chi$ of a 
quadratic extension $M$ 
of $K$, which we denote by ${\rm Ind}_M^K(\chi)$.) Anyhow we manage to solve 
this problem 
and establish the following result, where
Res$_M^K$ denotes, for any extension $M/K$ of number fields, 
the restriction of representations of $W_K$ to $W_M$: 

\medskip

\noindent{\bf Theorem D} \, \it Let $K/F$ be a quadratic extension of number 
fields with non-trivial automorphism 
$\theta$, and let $\pi$ be a 
cuspidal automorphic representation of GL$(2, \A_K)$. Then
\begin{itemize}
\item[(a)]{$As(\pi)$ is automorphic;}
\item[(b)]{If $\pi$ is non-dihedral, then $As(\pi)$ is a cusp form iff $\pi \circ \theta$ is not 
isomorphic to $\pi \otimes \mu$, 
for any idele class character $\mu$ of $K$;}
\item[(c)]{If $\pi$ is dihedral, i.e., associated to
a representation $\sigma = {\rm Ind}_M^K(\chi)$ of $W_K$ for a 
character $\chi$ of a quadratic extension $M$ of $K$, then $As(\pi)$ 
is cuspidal iff $M/F$ is non-Galois and
the representation 
${\rm Res}_M^K(\sigma^\theta) \otimes \chi$ 
of $W_M$ does not extend to a 
representation of $W_F$.}
\end{itemize}
\rm

\medskip

It may be helpful for the reader to note the following concrete description of the {\it Asai $L$-function} in a special case. Suppose $F = \Q$, $K$ a real quadratic field of class number $1$ with ring of integers $\mathfrak O_K$, and $\pi$  the representation defined by a holomorphic Hilbert modular newform $f$ of weight $2$ with 
$$
L(s, \pi) \, = \, \sum_{\mathfrak a}{}^{\prime} \, c(\mathfrak a) N({\mathfrak a})^{-s}.
$$
Then we have
$$
L(s, As(\pi)) \, = \, \zeta(2s-2) \sum_{m \geq 1} c(m{\mathfrak O}_K) m^{-s}.
$$
The interest in this comes, on the analytic side, from the fact that one sums over a sparse subset of the non-zero  ideals $\mathfrak a$ of $\mathfrak O_K$ to get this $L$-function, and on the geometric side, from the fact that $L(s, As(\pi))$ is a factor of the degree $2$ $L$-function of the associated Hilbert modular surface (cf. [Ra2], sec. 4, for example). It will be natural if one is reminded of the symmetric square $L$-function of a cusp form on GL$(2)/\Q$, where effectively one sums over the squares of positive integers. Indeed, for any quadratic extension $K/F$ of number fields and for any cusp form $\pi$ on GL$(2)/K$, if $\pi$ is the base change of a form $\pi_0$ on GL$(2)/F$, $L(s, As(\pi))$ factors as $L(s, {\rm sym}^2(\pi_0) \otimes \delta)L(s, \omega_0\delta)$, where $\omega_0$ is the central character of $\pi_0$.

\medskip

We will now make some remarks about the proofs of these results.

In section 3 we show how to reduce the proof of Theorem A to the following

\medskip

\noindent{\bf Theorem A$^\prime$} \, \it Fix a quadratic extension $K/F$ of number fields with associated quadratic character $\delta$ of Gal$(\overline F/F)$. Let $\rho$ be an irreducible $4$-dimensional, solvable $\C$-representation of Gal$(\overline F/F)$ whose restriction $\rho_K$ to Gal$(\overline F/K)$ is a tensor product of two $2$-dimensional representations. Then 
\begin{enumerate}
\item[(a)]{ \, $\rho$ is modular, i.e., there is a cuspidal automorphic representation $\Pi$ of GL$(4, \A_F)$ such that
$L(s, \rho) = L(s, \Pi_f)$;} 
\item[(b)]{ \, If $\rho_K$ is irreducible, one of the following happens:
\begin{enumerate}
\item[(i)] \, $\rho \, \simeq \, \tau \otimes \tau'$ over $F$, with dim$(\tau) = $dim$(\tau') = 2$;
\item[(ii)] \, $\rho \, \simeq \, $Ind$_L^F(\eta)$, with $L/F$ quadratic, $L \ne K$, and $\eta$ a $2$-dimensional representation of Gal$(\overline F/L)$; 
\item[(iii)] \, $\rho \, \simeq \, As(\sigma) \otimes \beta$, with $\sigma$ a $2$-dimensional representation of Gal$(\overline F/K)$ and $\beta$ a character of Gal$(\overline F/F)$.
\end{enumerate}}
\end{enumerate}
\rm

\medskip

Then we show, still in section 3, how to deduce Theorem A$^\prime$ modulo Theorem D; in fact we need Theorem D only in the (crucial) case when (iii) occurs. When (i) or (ii) occurs, the modularity already follows from Theorem M of [Ra1] and base change ([AC]), with the desired $\Pi$ being in case (i) (resp. (ii)) the Rankin-Selberg product $\pi(\tau) \otimes \pi(\tau')$ (resp. the automorphic induction $I_L^F(\pi(\eta))$); here $\tau \to \pi(\tau)$ the Langlands-Tunnell map on solvable $2$-dimensional Galois representations.

\medskip

The proof of Theorem D is accomplished in sections 4 through 7. The approach is similar to, but somewhat more subtle than, the proof of the 
existence of $\boxtimes: \, {\mathcal A}($GL$(2)/F) \times {\mathcal A}($GL$(2)/F) \to 
{\mathcal A}($GL$(4)/F)$ in [Ra1]. (For any $n \geq 1$, ${\mathcal A}($GL$(n)/F)$ denotes the set of isomorphism
classes of isobaric automorphic representations of GL$(n, \A_F)$.) 
In section 4 we show why it suffices to have the requisite properties at almost all places. Then in section 5 the {\it distinguished} case, i.e., when $\pi \circ \theta$ is an abelian twist of
$\pi$, is treated separately. In the {\it general} situation, i.e., when 
$\pi \boxtimes (\pi \circ \theta)$ is cuspidal,
we crucially use the converse 
theorem for GL$(4)$ due to Cogdell and Piatetski-Shapiro ([CoPS1], which requires knowledge of the 
{\it niceness} of the
twisted $L$-functions $L(s, As(\pi) \times \pi')$ for automorphic forms $\pi'$ of GL$(2)/F$ for a 
suitable class of $\pi'$. 
Many properties of 
certain closely related functions, to be denoted $L_1(s, As(\pi) \times \pi')$, were established by 
Piatetski-Shapiro and Rallis ([PS-R]), and by Ikeda 
([Ik1,2]), via an integral representation, which we use. There is another possible approach to 
studying $L(s, As(\pi) \times \pi')$ via the Langlands-Shahidi method [Sh1], which yields another 
family of closely related $L$-functions, denoted $L_2(s, As(\pi) \times \pi')$, 
with boundedness properties established recently in [GeSh], 
but we will not use it and our argument here is hewed to make 
use of the integral representation. 

For almost all finite places $v$, the local factors 
$L(s, As(\pi_v) \times \pi_v')$ and $L_1(s, As(\pi_v) \times \pi_v')$ agree.
But a thorny problem arises however, due to our inability to identify the bad and archimedean factors. 
In fact, when $F_v = \R$, one does not even have a computation of the corresponding $L_1$-factor when 
$\pi, \pi'$ are unramified.
Luckily, things simplify quite a bit under suitable, solvable base changes $K/F$ with $K$ {\it totally 
complex}, and after 
constructing the base-changed candidates $As(\pi)_K$
for an infinitude of such $K$, we descend to $F$ as in sections 3.6, 3.7 of [Ra1].  
We also have to control the intersection of the 
ramification loci of $As(\pi)$, $\pi'$ and $K/F$. 

One of the reasons why we have to work with $L(s, As(\pi) \times \pi')$ is that we know how its local 
$\varepsilon$-factors behave under twisting by a highly ramified character, 
and this is not apriori the case with 
$L_j(s, As(\pi) \times \pi')$ for $j = 1$ or $2$. Indeed this will present a difficulty for the lifting of
(generic) automorphic forms from GO$(2n)$ to GL$(2n)$, and for large $n$ one cannot, as of yet, 
make use of base change and descent as we do here. For the lifting from GO$(2n+1)$ to GL$(2n)$, see [CoKPSSh].

\vskip 0.08in

We would like to thank R.P. Langlands and the Institute for Advanced Study, Princeton, for their 
hospitality during the year 1999-2000, and the American Institute of 
Mathematics, Palo Alto, for support during the month of August in 1999. 
This project was partially funded by the NSF and AMIAS. We 
also want to acknowledge with thanks
some helpful conversations with Michael Aschbacher and David Wales
concerning the representations of finite groups. The main result here was presented at a 
conference on Automorphic Forms in Luminy, France (May 1999), and then at the IAS, Princeton (April 2000),
and ICTP, Trieste, Italy (August 2000), where we received some helpful feedback.
We would like to thank a number of people, 
including Jim Arthur, Don Blasius, Kevin Buzzard, Jim Cogdell, Bill Duke, Farshid Hajir, Henry Kim, Jeff Lagarias, Nick Katz, Dipendra Prasad, A.Raghuram, Peter Sarnak, Jean-Pierre Serre, Freydoon Shahidi and Joe Shalika who have shown interest and/or made comments.

\vskip 0.2in

\section{\bf Preliminaries on orthogonal similitude groups}

\bigskip

Here we collect some basic facts, which we will need.

\medskip 

Let $k$ be a field of characteristic different from $2$, with separable algebraic closure $\overline k$.
If $V$ is a finite dimensional vector space with a non-degenerate, symmetric bilinear form $B$, the 
associated orthogonal similitude group is
$$
{\text GO}(V,B) \, := \, \{g \in {\text GL}(V) \, \vert \, B(gv, gw) \, = \, \lambda(g)B(v,w), \,
{\text with} \, \lambda(g) \in k^\ast, \, 
\forall v,w \in V\}.
\leqno(1.1)
$$
The character $\lambda: {\text GO}(V,B) \, \rightarrow \, k^\ast, \, g \to \lambda(g)$, 
is the {\it similitude factor}. The kernel of $\lambda$ is the {\it orthogonal group} O$(V,B)$,
whose elements necessarily have determinant $\pm 1$, and the kernel of det is the
{\it special orthogonal group} SO$(V,B)$. 

If $V = k^n$ with $B$ the standard bilinear form $B_0: (v,w) \to {}^tvw$, then one writes GO$(n,k)$,
O$(n,k)$ and SO$(n,k)$ instead of
GO$(V,B)$, O$(V,B)$ and SO$(V,B)$. 
Denote by $Z_n(k)$ the {\it center} of GO$(n,k)$ consisting of all the scalar matrices $cI_n$, $c \in k^\ast$.
Clearly, $\lambda(cI_n) = c^2$, so that ${k^\ast}^2$ is in the image of $\lambda$.
The {\it odd dimensional case} is relatively simple. One has

\medskip

\noindent{\bf Lemma 1.2} \, \it If $n$ is odd and $k = \overline k$, then we have the direct product
decomposition
$$
{\text GO}(n,k) \, = \, {\text SO}(n,k) \times Z_n(k).
$$
\rm

\medskip

Indeed, as $k = \overline k$, $\lambda(Z_n(k))$ is all of $k^\ast$, and since O$(n,k)$ is by definition 
the kernel of $\lambda$, GO$(n,k)$ is generated by the normal subgroups
O$(n,k)$ and $Z_n(k)$. On the other hand, the 
intersection of these two groups is simply $\{\pm I_n\}$. Since $n$ is odd, the image of
${\sl det}: {\text O}(n,k) \, \rightarrow \, \{\pm 1\}$ is the same as that of $\{\pm I_n\}$. 
The assertion follows.

\medskip

Note that GO$(1,k) = Z_1(k) = k^\ast$. There is a useful description in the $n=3$ case, which we will
now recall. The {\it adjoint representation}
$$
{\text Ad}: \, {\text PGL}(2, k) \, \rightarrow \, {\text GL}(3, k)
$$
is irreducible and self-dual with determinant $1$. This identifies the image of Ad with SO$(3,k)$,
thanks to the simplicity of the latter. 
By abuse 
of notation, we will also write Ad for its composition with the canonical map from GL$(2,k)$ 
onto PGL$(2,k)$. This gives rise to the short exact sequence
$$
1 \rightarrow k^\ast \rightarrow {\text GL}(2,k) \rightarrow {\text SO}(3,k) \rightarrow 1,
\leqno(1.3)
$$
where the maps in the middle are $c \to cI_2$ and $g \to {\text Ad}(g)$.

\medskip

The {\it even dimensional case} $n = 2m$ is more interesting. 
Since for any $g$ in GO$(2m,k)$, the square of its determinant
is $\lambda(g)^{2m}$, we can define a homomorphism, called the 
{\it similitude norm}
$$
\nu: \, {\text GO}(2m,k) \, \rightarrow \, \{\pm 1\},
\leqno(1.4)
$$
by sending $g$ to $\lambda(g)^{-m}{\text det}(g)$. 

The kernel of $\nu$, denoted SGO$(2m,k)$, is called the 
{\it special orthogonal similitude group}. (Some people write GSO$(2m,k)$
instead.) The map $\nu$ does not split. 
Since $\nu$ is just the determinant map on O$(2m,k)$, the intersection
of SGO$(2m,k)$ with O$(2m,k)$ is SO$(2m,k)$.
When $k = \C$, SGO$(2m,k)$ (resp. SO$(2m,k)$) is
the connected component of GO$(2m,k)$ (resp. O$(2m,k)$). SGO$(2,k)$ is a torus.

Note that $\nu(cI_{2m})$ is $1$, and that
SGO$(2m,k)$ is generated by SO$(2m,k)$ and $Z_{2m}(k)$;
but their intersection is $\{\pm I_{2m}\}$. 

\medskip

We will conclude this section by recalling a low dimesional isomorphism for 
$k = \overline k$, which we 
will need, between SGO$(4,k)$ and a quotient of GL$(2,k) \times
{\text GL}(2,k)$. 

Let $W$ be $k^2$ with the standard symplectic form given by the determinant. Then
the induced bilinear form $B$ on the tensor product $W \otimes W$ is non-degenerate and
symmetric. There is an isometry between $(W \otimes W, B)$ and $(k^4, B_0)$. Since
GL$(2,k)$ is the symplectic similitude group of $(W, {\text det})$, we
get an exact sequence
$$
1 \rightarrow k^\ast \rightarrow {\text GL}(2,k) \times {\text GL}(2,k)
\rightarrow GO(4,k),
\leqno(1.5)
$$
where the map on $k^\ast$ is just given by $c \to (cI_2, c^{-1}I_2)$. The map $\beta$, say, 
on the right can be 
described explicitly as follows. The quadratic space $(k^4, B_0)$ is also isometric to
$(M_2(k), B_1)$, where $B_1$ is the symmetric bilinear map $(X,Y) \to {\rm tr}({}^tXY)$. Under this
identification, $\beta(g,g')$ is, for all $g, g'$ in GL$(2,k)$, the automorphism of
$k^4$ given by $X \to {}^tgXg'$. Clearly the kernel of $\beta$ consists of pairs
$(cI_2, c^{-1}I_2)$ with $c \in k^\ast$, proving the requisite exactness. 

Note that $\lambda(\beta(g,g'))$ is det$(g)$det$(g')$, while the determinant of
$\beta(g,g')$ is its square. Hence $\nu$ is trivial on the image of $\beta$. It is easy
to see that $Z_4(k)$ lies in the image of $\beta$, and that 
$\beta({\rm SL}(2,k) \times {\rm SL}(2, k))$ is a subgroup of SO$(4,k)$ properly containing
$\{\pm I_4\}$. The abelianization of SO$(4,k)$ is $k^\ast/{k^\ast}^2$ ([D], p. 57), and since 
we have assumed that $k = \overline k$, SO$(4, k)$ is perfect, i.e., it equals its own commutator 
subgroup. Then by the discussion on page 
59 of {\it loc. cit}, SO$(4, k)/\{\pm I_4\}$ is isomorphic to PSL$(2, k) \times {\rm PSL}(2,k)$.
It follows that $\beta$ maps SL$(2, k) \times $SL$(2, k)$ onto SO$(4, k)$ with kernel 
$\{\pm(I_2,I_2)\}$. Putting everything together, we obtain
$$
\beta({\text GL}(2,k) \times {\text GL}(2,k)) \, = \, {\text SGO}(4,k).
\leqno(1.6)
$$

\vskip 0.2in

\section{\bf The reducible case}

\bigskip

Suppose we are given a representation $\rho$ as in the statement of Theorem A,
which is reducible. Thanks to Maschke's theorem we may
write $\rho \simeq \oplus_j \rho_j$, with each
$\rho_j$ irreducible of dimension $n_j$, and $\sum_j n_j = 4$.
Suppose we have found, for each $j$, a cuspidal automorphic representation
$\pi_j = \pi_{j, \infty} \otimes \pi_{j,f}$ of GL$(n_j, \A_F)$ 
such that $L(s, \rho_j) = L(s, \pi_{j,f})$. Then we can 
consider the {\it isobaric sum} of Langlands ([La2], [JS])
$$
\pi \, = \, \boxplus_j \pi_j,
\leqno(2.1)
$$
which is automorphic and satisfies 
$$
L(s, \pi) \, = \, \prod_j L(s, \pi_j).
$$
Since the $L$-functions of Artin are also additive, we get
$L(s, \rho) = L(s, \pi_f)$ as desired. So it remains to find the $\pi_j$.

Note that cuspidal automorphic representations of GL$(1,\A_F)$ are just idele class characters
of $F$. So when $n_j = 1$, the existence of $\pi_j$ follows from class field theory.

Since the image of $\rho$ is by hypothesis solvable, the same will be true for each 
$\rho_j$. So if $n_j = 2$, we may apply the celebrated theorem of Langlands ([La1]) and Tunnell
([Tu]) to conclude the existence of $\pi_j$. 

It remains to consider the case when $n_j$ is $3$ for some $j$, say for $j=1$. Then we must
have a decomposition
$$
\rho \, \simeq \, \rho_1 \oplus \rho_2,
$$
with $\rho_1$ (resp. $\rho_2$) irreducible of dimension $3$ (resp. $1$). Since by hypothesis,
the image of $\rho$ lands is GO$(4, \C)$, and since there can be no intertwining 
between $\rho_1$ and $\rho_j$, we must have
$$
im(\rho_1) \, \subset \, {\text GO}(3, \C).
$$
Thanks to Lemma 1.2, GO$(3, \C)$ is SO$(3, \C) \times \C^\ast$. So we may
write
$$
\rho_1 \, \simeq \, \rho' \otimes \chi,
\leqno(2.2)
$$
where $\chi$ is a character ${\mathfrak G}_F \rightarrow \C^\ast$, and $\rho'$ 
is a $3$-dimensional representation of ${\mathfrak G}_F$ with image in
SO$(3, \C)$. 

Moreover,
the exact sequence (1.3), which can be viewed as an exact sequence of
trivial modules under ${\mathfrak G}_F =$ Gal$(\overline F/F)$, 
furnishes the cohomology exact sequence
$$
{\text Hom}({\mathfrak G}_F, {\text GL}(2, \C)) \rightarrow 
{\text Hom}({\mathfrak G}_F, {\text SO}(3, \C)) \rightarrow
H^2({\mathfrak G}_F, \C^\ast),
\leqno(2.3)
$$
with $\rho'$ belonging to the middle group. On the other hand, a theorem of Tate
(see [Se], for a proof) asserts that the group on the right hand side of (2.3) is
trivial as $F$ is a number field. Thus we may lift $\rho'$ to an element of the left hand
side group of (2.3). In other words, we can find a (non-unique) $2$-dimensional representation
$\tau_1$ of ${\mathfrak G}_F$ such that 
$$
\rho_1 \, \simeq \, {\text Ad}(\tau_1) \otimes \chi.
\leqno(2.4)
$$
Since $\rho_1$ has solvable image, $\tau_1$ is also forced to have the same property. Applying
Langlands-Tunnell once again, we get an isobaric automorphic representation $\eta_1$, which 
must in fact be cuspidal as $\rho_1$ and hence $\tau_1$ are irreducible, such that
$L(s, \tau_1)$ equals $L(s, \eta_{1,f})$. 

By [GeJ] one knows that, given any cuspidal automorphic representation 
$\eta$ of GL$(2, \A_F)$, there exists
a functorially associated (isobaric) automorphic representation Ad$(\eta)$ such that
$$
L(s, {\text Ad}(\eta)) \, = \, \prod_v L(s, {\text Ad}(\sigma_v(\eta))),
\leqno(2.5)
$$
where the product is over all the places $v$ of $F$, and $\sigma_v(\eta)$ 
(resp. $\sigma(\eta_v)$) is the 
$3$-dimensional representation of the Weil group (resp. Weil-Deligne group) $W_{F_v}$
(resp. $W'_{F_v}$) when $v$ is archimedean (resp. non-archimedean), 
associated to $\eta_v$ by the local Langlands correspondence for GL$(n)$ ([HaT], [He]).

Then it follows that 
$$
L(s, \rho_1) = L(s, ({\text Ad}(\eta_1) \otimes \chi)_f).
$$
So we are done by setting $\pi_1 = {\text Ad}(\eta_1) \otimes \chi$. 

\qed

\vskip 0.2in

\section
{\bf Modularity modulo Theorem D}

\bigskip

In this section we will show how to prove Theorem A if we admit the truth of Theorem D. 
We will also need to make use of the main theorem of [Ra1].
Thanks to the discussion in the previous section, we may assume that $\rho$ is irreducible.

\medskip

By hypothesis, the image of $\rho$ lies in GO$(4, \C)$. Recall from section 1 the definion of the 
subgroup SGO$(4, \C)$, which is the kernel of the {\it similitude norm}
$$
\nu: \, {\text GO}(4, \C) \, \rightarrow \, \{\pm 1\}.
$$

Let $K$ be the extension of $F$ defined by the kernel of $\nu \circ \rho$. 
Then $[K : F] \leq 2$.
Write $\rho_K$ for the restriction of $\rho$ to 
${\mathfrak G}_K =$ Gal$(\overline F/K)$. Thanks to (1.5) and (1.6),
one has the following short equence
of trivial Galois modules:
$$
1 \rightarrow \C^\ast \rightarrow {\text GL}(2,\C) \times {\text GL}(2,\C)
\rightarrow {\text SGO}(4,\C) \rightarrow 1,
\leqno(3.1)
$$
where the maps in the middle are $c \to (cI_2, c^{-1}I_2)$ and $(g,g')
\to (X \to {}^tgXg')$. The associated (continuous)
cohomology exact sequence gives
$$
{\text Hom}({\mathfrak G}_K, \C^\ast) \rightarrow
{\text Hom}({\mathfrak G}_K, {\text GL}(2, \C) \times {\text GL}(2, \C)) \rightarrow 
{\text Hom}({\mathfrak G}_K, {\text SGO}(4, \C)) \rightarrow
H^2({\mathfrak G}_K, \C^\ast),
\leqno(3.2)
$$
with $\rho_K$ belonging to the second group from the right. Recall Tate's theorem which says
that the first group on the right is trivial. So we may find $\sigma, \sigma'$
in ${\text Hom}({\mathfrak G}_K, {\text GL}(2, \C))$ such that
$$
\rho \, \simeq \, \sigma \otimes \sigma'.
\leqno(3.3)
$$
Note that this lifting is unique only up
to flipping the two factors and
changing $(\sigma, \sigma')$ by $(\sigma \otimes \mu, \sigma' \otimes \mu^{-1})$, for any
character $\mu \in {\text Hom}({\mathfrak G}_K, \C^\ast)$. 

\medskip

Since the image of $\rho$ was assumed to be solvable, we see easily that the images of $\sigma, \sigma'$
should also be solvable. And since $\rho$ is irreducible, the same should hold for
$\sigma$ and $\sigma'$. So we may apply the theorem of Langlands and Tunnell to deduce the existence
of cuspidal automorphic representations $\pi, \pi'$ of GL$(2, \A_F)$, 
respectively associated to $\sigma, \sigma'$.

\medskip

Now suppose the image of $\rho$ lands in SGO$(4, \C)$ itself, in which case $K = F$. Then by 
Theorem M (in section 3)
of [Ra1], we know the existence of an isobaric automorphic representation $\pi \boxtimes \pi'$
of GL$(4, \A_F)$ such that 
$$
L(s, \pi \boxtimes \pi') \, = \, L(s, \pi \times \pi'),
\leqno(3.4)
$$
where the $L$-function on the right is the Rankin-Selberg $L$-function associated to the pair $(\pi, \pi')$.
In addition, we have at any place $v$,
the local factors of $\pi \boxtimes \pi'$ identify functorially with those of
the tensor product $\sigma_v(\pi) \otimes \sigma_v(\pi')$ of the local Langlands parameters
$\sigma(\pi), \sigma(\pi')$ ([HaT], [He]), proved long ago for GL$(2)$ by P. Kutzko. Since $\pi$ 
(resp. $\pi'$) is associated to $\sigma$ (resp. 
$\sigma'$), the local representations $\sigma_v(\pi)$ (resp. $\sigma_v(\pi')$) are isomorphic to
the ones defined by the restriction at $v$ of $\sigma$ (resp. $\sigma'$). Thus the
automorphic representation $\Pi$ of GL$(4, \A_F)$, whose existence is predicted by Theorem A, is none
other than $\pi \boxtimes \pi'$. The cuspidality criterion of [Ra1] (Theorem M) shows easily that, since 
$\rho$ is irreducible, $\Pi$ must be cuspidal. 
We are done in this case.

\medskip

We may henceforth assume that $[K:F] = 2$, which is the subtler case. Denote by $\theta$ the non-trivial
automorphism of $K$ over $F$. We can again find cuspidal 
automorphic representations $\pi, \pi'$ of GL$(2, \A_K)$ such that
$$
L(s, \rho_K) \, = \, L(s, \pi_f \boxtimes \pi'_f),
\leqno(3.5)
$$
which proves that the restriction $\rho_K$ of $\rho$
to ${\mathfrak G}_K$ is modular. 

\medskip

We will now explain why this case is difficult.
The identity (3.5) implies that $\pi \boxtimes \pi'$ is $\theta$-invariant, so by the base change theorem of 
Arthur and Clozel ([AC]), we can find an isobaric automorphic representation $\Pi$ of
GL$(4, \A_F)$ such that its base change $\Pi_K$ is isomorphic to
$\pi \boxtimes \pi'$. One can also see easily that the local factors of $\Pi$ and $\rho$ agree at
all the places of $F$ which split in $K$. But one is stuck at this point and cannot easily deduce the
requisite identity at the {\it inert} places, except when $\rho_K$ is no longer irreducible.

\medskip

It is now clear that Theorem A is a consequence of Theorem $A\prime$. So we will address the following:

\medskip

{\it Proof of Theorem A$^\prime$ modulo Theorem D}:

Suppose $\rho_K$ is {\it reducible}. Then it must contain an irreducible summand $\tau$ of 
dimension $\leq 2$. By
Frobenius reciprocity, $\rho$ should intertwine with the induction ${\text Ind}_K^F(\tau)$ of $\tau$
to ${\mathfrak G}_F$. As $\rho$ is irreducible of dimension $4$ and ${\text Ind}_K^F(\tau)$ is at most 
of dimension $4$, we are forced to have
$$
\rho \, \simeq \, {\text Ind}_K^F(\tau),
$$
with dim$(\tau) = 2$. The solvability of the image of $\rho$ implies the same about
that of $\tau$, and so we may apply Langlands-Tunnell to get a cuspidal automorphic
representation $\eta$ of GL$(2, \A_K)$ associated to $\tau$, and we are done by taking $\Pi$ to be the
automorphically induced representation $I_K^F(\tau)$ (see [AC], and also [Ra1], sec. 2). Since $\rho$ is irreducible, $\tau$, and hence $\eta$, cannot be $\theta$-invariant, where $\theta$ is the non-trivial automorphism of $K/F$. Consequently, $I_K^F(\eta)$ must be cuspidal.

\medskip

So we may, and we will, assume that $\rho_K$ is irreducible. 
Since it is the restriction of $\rho$,
we have
$$
(\sigma \otimes \sigma')^{\theta} \, \simeq \, \rho_K^{\theta} \, \simeq \, \rho_K 
\, \simeq \sigma \otimes \sigma^{\theta},
$$
where $\theta$ is the nontrivial automorphism of $K$ over $F$. 

\medskip

We will first prove part (b) of Theorem A$^\prime$.

\medskip

\noindent{\bf Lemma 3.6} \, The irreducibility of $\rho_K = \sigma \otimes \sigma'$ implies
that at least one of the representations $\sigma, \sigma'$ is non-dihedral, and $\sigma'$ cannot be a one-dimensional twist of $\sigma$. Furthermore, since $\rho_K$ is $\theta$-invariant, one of the following must happen:
\begin{enumerate}
\item[(I)] \, There exists a character $\nu$ of Gal$(\overline F/K)$ such that
$$
\sigma^\theta \, \simeq \, \sigma \otimes \nu \quad {\rm and} \quad 
(\sigma')^\theta \, \simeq \, \sigma' \otimes \nu^{-1};
$$
\item[(II)] \, There exists a character $\mu$ of Gal$(\overline F/K)$ such that
$$
\sigma^{\theta} \, \simeq \, \sigma' \otimes \mu \quad {\rm and} \quad
(\sigma')^{\theta} \simeq \sigma \otimes \mu^{-1}.
$$
\end{enumerate}
\rm

\medskip

{\it Proof of Lemma 3.6}. \, Let $\omega$, resp. $\omega'$, denote the determinant of $\sigma$, resp. $\sigma'$. Clearly, if $\sigma' \simeq \sigma \otimes \lambda$ for some character $\lambda$, then $\rho_K$ is reducible, as $\sigma \otimes \sigma'$ will then be $({\rm sym}^2(\sigma) \otimes \lambda) \oplus \omega\lambda$. We will now show that not both $\sigma, \sigma'$ can be dihedral.
We have
$$
{\rm End}(\rho_K) \, \simeq \, \rho_K \otimes \rho_K^\vee \, \simeq \, \rho_K \otimes (\sigma \otimes \omega^{-1}) \otimes (\sigma' \otimes {\omega'}^{-1}) \, \simeq \, \rho_K^{\otimes 2} \otimes (\omega\omega')^{-1}.
$$
Since $\rho_K$ is irreducible, any linear character occurring in End$(\rho_K)$ must have multiplicity one. We claim that if $\sigma$ and $\sigma'$ are both dihedral, then End$(\rho_K)$ contains the trivial representation with multiplicity 
$> 1$. Indeed suppose $\sigma = {\rm Ind}_M^K(\chi)$ and $\sigma' = {\rm Ind}_N^K(\chi')$, where $M, N$ are quadratic extensions of $K$, and $\chi, \chi'$ characters of $\mathfrak G_M$, $\mathfrak G_N$ respectively. Since $\rho_K^{\otimes 2}$ is sym$^2(\rho_K) \oplus \Lambda^2(\rho_K)$, it suffices to prove that $\omega\omega'$
occurs twice in sym$^2(\rho_K)$. Note that 
$$
{\rm sym}^2(\sigma) \, \simeq \, {\rm Ind}_M^K(\chi^2) \oplus \omega\epsilon
\quad {\rm and} \quad 
{\rm sym}^2(\sigma') \, \simeq \, {\rm Ind}_M^K({\chi'}^2) \oplus \omega'\epsilon',
$$
where $\epsilon$ (resp. $\epsilon'$) is the quadratic character of $\mathfrak G_K$ corresponding to $M/K$ (resp. $N/K$). Since
$$
{\rm sym}^2(\sigma \otimes \sigma') \, \simeq \, {\rm sym}^2(\sigma) \otimes {\rm sym}^2(\sigma') \oplus \Lambda^2(\sigma) \otimes \Lambda^2(\sigma'),
\leqno(3.7)
$$
we get
$$
{\rm sym}^2(\rho_K) \, \simeq \, \simeq \, {\rm Ind}_M^K(\chi^2) \otimes {\rm Ind}_M^K({\chi'}^2)
\oplus {\rm Ind}_M^K(\chi^2)\otimes \omega'\epsilon' \oplus {\rm Ind}_M^K({\chi'}^2)\otimes \omega\epsilon \oplus (\omega\omega')^{\oplus 2},
$$
as asserted. So at least one of $\sigma, \sigma'$ must be non-dihedral.
Finally by $\theta$-invariance,
$$
\sigma^\theta \otimes (\sigma')^\theta \, \simeq \, \sigma \otimes \sigma',
$$
which evidently implies, by using irreducibility, that we must be in case (I) or (II).

\qed

\medskip

\noindent{\bf Lemma 3.8} \, \it Suppose we are in case (I) of Lemma 3.6. Then there exist characters $\chi, \chi'$ of $\mathfrak G_K$ such that 
$\sigma \otimes \chi$ and $\sigma' \otimes \chi'$ are both $\theta$-invariant.
\rm

\medskip

{\it Proof of Lemma 3.8}. \, We may assume, after interchanging $\sigma, \sigma'$ if necessary, that $\sigma'$ is non-dihedral. We claim that  
$$
(\omega\omega')^\theta \, = \, \omega\omega'.
\leqno(3.9)
$$
To begin, note that the $\theta$-invariance of $\rho_K$ implies the same for its contragredient and its symmetric and exterior powers. The determinant of $\rho_K$ is easily seen to be $(\omega\omega')^2$, and this shows that
$$
\beta: \, = \, (\omega\omega')^\theta/\omega\omega'
$$
has order $\leq 2$. Suppose the order is $2$. Then the identity (3.7) shows that 
$$
\beta \, \subset \, {\rm Ad}(\sigma) \otimes {\rm Ad}(\sigma'),
$$
where for any $2$-dimensional $\tau$, 
$$
{\rm Ad}(\tau): \, = \, {\rm sym}^2(\tau) \otimes \omega_\tau^{-1},
$$
the adjoint representation.
The self-duality of Ad$(\sigma)$, plus the irreducibility of Ad$(\sigma')$, then implies that Ad$(\sigma')$ is isomorphic to Ad$(\sigma) \otimes \beta$, which is the same, as $\beta^2 = 1$, as Ad$(\sigma \otimes \beta)$.
But we have the following

\medskip

{\bf Theorem 3.10} \, \it Let $\tau, \tau'$ be irreducible, $2$-dimensional representations of $\mathfrak G_K$ with isomorphic adjoint representations. Then there exists a character $\chi$ of $\mathfrak G_K$ such that
$$
\tau' \, \simeq \tau \otimes \chi.
$$
\rm

\medskip

For a proof see [Ra2]; it is important to note that it holds whether or not $\tau$ is dihedral. This is the Galois version of the so called {\it multiplicity one for SL$(2)$}. Its automorphic version was proved in [Ra1] (see Theorem 4.1.2).

Applying this with $\tau = \sigma \otimes \beta$ and $\tau' = \sigma'$, we see that $\sigma'$ is a one dimensional twist of $\sigma$, which contradicts the irreducibility of $\rho_K$. Hence $\beta$ must be $1$ and the claim is proved.

\medskip

Next we claim that
$$
({\rm sym}^2(\sigma) \otimes \omega')^\theta \, \simeq \, {\rm sym}^2(\sigma) \otimes \omega' \quad
{\rm and} \quad
({\rm sym}^2(\sigma') \otimes \omega)^\theta \, \simeq \, {\rm sym}^2(\sigma') \otimes \omega.
\leqno(3.11)
$$
Indeed, since $\sigma^\theta$ is by hypothesis $\sigma \otimes \mu$, we have
$$
({\rm sym}^2(\sigma) \otimes \omega')^\theta \, \simeq \, {\rm sym}^2(\sigma \otimes \mu) \otimes (\omega')^\theta \, \simeq \, ({\rm sym}^2(\sigma) \otimes \omega')\otimes \mu^2\left(\frac{(\omega')^\theta}{\omega'}\right).
\leqno(3.12)
$$
On the other hand, comparing the determinants of $\sigma^\theta$ and $\sigma \otimes \mu$, we get
$$
\mu^2 \, = \, \frac{\omega^\theta}{\omega}.
\leqno(3.13)
$$
The first half of the asserted identity (3.11) now follows by combining (3.9), (3.12) and (3.13). The proof of the second half is the same.

\medskip

We will now prove the $\theta$-invariance of $\sigma \otimes \chi$ for a suitable $\chi$. The case of $\sigma'$ is similar and will be left to the reader.

Note that sym$^2(\sigma \otimes \omega')$ is of GO$(3)$-type. Explicitly,
$$
{\rm sym}^2({\rm sym}^2(\sigma) \otimes \omega') \, \simeq \, ({\rm sym}^4(\sigma) \oplus \omega^2) \otimes {\omega'}^2,
$$
and so $(\omega\omega')^2$ occurs in the symmetric square of sym$^2(\sigma) \otimes \omega'$. Since $\omega\omega'$ is $\theta$-invariant by (3.9), it extends to a character of $\mathfrak G_F$. Moreover, the $\theta$-invariance of sym$^2(\sigma) \otimes \omega'$ (cf. (3.11)) gives the existence of a $3$-dimensional representation $\eta$ of $\mathfrak G_F$ such that
$$
{\rm Res}_K^F(\eta) \, \simeq \, {\rm sym}^2(\sigma) \otimes \omega'.
$$
Since the restriction of the symmetric square of $\eta$ to $\mathfrak G_K$ contains the $\theta$-invariant character $(\omega\omega')^2$, there will be a character $\lambda$, say, of $\mathfrak G_F$, such that
$$
\lambda \, \subset \, {\rm sym}^2(\eta) \quad {\rm and} \quad {\rm Res}_K^F(\lambda) = (\omega\omega')^2.
$$
In other words, $\eta$ is also of GO$(3)$-type. By Lemma 1.2 and the exact sequence (1.3), we then get the existence of a $2$-dimensional representation $\tau$, and a character $\alpha$, of $\mathfrak G_F$ such that
$$
\eta \, \simeq \, {\rm Ad}(\tau) \otimes \alpha.
$$
This yields the isomorphism
$$
{\rm Ad}({\rm Res}_K^F(\tau)) \otimes {\rm Res}_K^F(\alpha) \, \simeq \, 
{\rm sym}^2(\sigma) \otimes \omega'.
$$ 
Since sym$^2(\sigma)$ is Ad$(\sigma) \otimes \omega$, we get
$$
{\rm Ad}(\sigma) \, \simeq \, {\rm Ad}({\rm Res}_K^F(\tau \otimes \nu)),
$$
for a character $\nu$ of $\mathfrak G_F$ satisfying
$$
{\rm Res}_K^F(\nu/\alpha) \, = \, (\omega\omega')^{-1}.
$$
The existence of $\nu$ comes from the $\theta$-invariance of $\omega\omega'$.
Applying Theorem 3.10 again, now with $\tau = \sigma$ and $\tau' =  {\rm Res}_K^F(\tau \otimes \nu)$, we get the existence of a character $\chi$ of $\mathfrak G_F$ such that
$$
\sigma \otimes \chi \, \simeq \, {\rm Res}_K^F(\tau \otimes \nu).
$$
Done with the proof of Lemma 3.8.

\qed

\medskip

\noindent{\bf Lemma 3.14} \, \it Suppose we are in case (I) of Lemma 3.6. Then the following hold:
\begin{enumerate}
\item[(i)] \, If $\sigma, \sigma'$ are both non-dihedral, then there exist irreducible, $2$-dimensional representations $\tau, \tau'$ of $\mathfrak G_F$ such that
$$
\rho \, \simeq \, \tau \otimes \tau';
\leqno(*)
$$
\item[(ii)] \, If either $\sigma$ or $\sigma'$ is dihedral, then either (*) holds or there exists a quadratic extension $L/F$ and a $2$-dimensional representation $\eta$ of $\mathfrak G_L$ such that
$$
\rho \, \simeq \, {\rm Ind}_L^F(\eta).
$$
\end{enumerate}
\rm

\medskip

{\it Proof of Lemma 3.14}. \, Again by the irreducibility of $\sigma \otimes \sigma'$, we may assume that $\sigma'$ is non-dihedral. By Lemma 3.8, there exist characters $\chi, \chi'$of $\mathfrak G_K$ such that $\sigma \otimes \chi$ and $\sigma' \otimes \chi'$ are $\theta$-invariant. We first claim that
$$
\nu^2 \, = \, 1 \quad {\rm where} \quad \nu \, = \, (\chi\chi')^\theta/\chi\chi'.
\leqno(3.15)
$$
Indeed, since we are in case (I), $\sigma^\theta$ is of the form $\sigma \otimes \mu$ for some character $\mu$, so that
$$
\sigma \, \simeq \, (\sigma \otimes \chi)^\theta \otimes \chi^{-1} \, \simeq \, \sigma \otimes \mu\chi^\theta/\chi.
$$
When $\sigma$ is non-dihedral, we must have $\mu = \chi/\chi^\theta$, and when $\sigma$ is dihedral, $\mu\chi^\theta/\chi$ must be trivial or quadratic. So
$$
(\mu\chi^\theta/\chi)^2 \, = \, 1
$$
in either case. Similarly, $(\mu^{-1}{\chi'}^\theta/\chi')^2$ is $1$. The claimed identity (3.15) is then a consequence.

This argument in fact shows that if $\sigma$ is {\it non-dihedral}, then $\nu = 1$. So we may choose a character $\gamma$ of $\mathfrak G_F$ such that
$$
\chi\chi' \, = \, {\rm Res}_K^F(\gamma).
$$
On the other hand, since $\sigma \otimes \chi$ and $\sigma' \otimes \chi'$ are $\theta$-invariant, there are irreducible $2$-dimensional repesentations $\tau, \tau'$ of $\mathfrak G_F$ such that 
$$
\sigma \otimes \chi \, \simeq \, {\rm Res}_K^F(\tau) \quad {\rm and} \quad \sigma' \otimes \chi' \, \simeq \, {\rm Res}_K^F(\tau').
$$
Putting all these together we get
$$
{\rm Res}_K^F(\rho) \, \simeq \, {\rm Res}_K^F(\tau \otimes \tau' \otimes \gamma).
$$
Also, $\tau$ cannot be a one-dimensional twist of $\tau'$ as it would make $\sigma \otimes \sigma'$ reducible. Then it follows that
$$
\rho \, \simeq \, \tau \otimes (\tau' \otimes \gamma'),
$$
where $\gamma'$ is either $\gamma$ or $\gamma\delta$. So we get (i).

So we may assume (for this Lemma) that $\sigma$ is {\it dihedral} and that $\nu$ is a non-trivial quadratic character of $\mathfrak G_K$. The $\theta$-invariance of $\sigma \otimes \sigma' \otimes \chi\chi'$ then implies that 
$\sigma \otimes \sigma'$ is isomorphic to $\sigma \otimes \sigma' \otimes \nu$. Consequently,
$$
\nu \, \subset \, (\sigma^\vee \otimes \sigma) \otimes ({\sigma'}^\vee \otimes \sigma') \, \simeq \, {\rm Ad}(\sigma) \otimes {\rm Ad}(\sigma') \oplus {\rm Ad}(\sigma) \oplus {\rm Ad}(\sigma') \oplus 1.
$$
Recall that for any irreducible $2$-dimensional representation $\tau$, Ad$(\tau)$ is reducible iff $\tau$ is dihedral. Since $\sigma'$ (resp. $\sigma$) is non-dhedral (resp. dihedral), we see that Ad$(\sigma) \otimes $Ad$(\sigma')$ and Ad$(\sigma')$ have no one-dimensional summands. So we must have
$$
\nu \, \subset \, {\rm Ad}(\sigma).
$$
This implies that $\sigma \otimes \nu \simeq \sigma$, i.e., that $\sigma$ is induced by a character of $\mathfrak G_M$, if $M$ denotes the quadratic extension of $K$ cut out by $\nu$. On the other hand, since $\theta^2 = 1$ and $\nu^2 = 1$,
$$
\nu^\theta \, = \, \left(\frac{(\chi\chi')^\theta}{\chi\chi'}\right) \, = \, \nu^{-1} \, = \, \nu.
$$
So $\nu$ is the restriction to $K$ of a character $\nu_0$, say, of $\mathfrak G_F$. Denote by $M_0$ the quadratic extension of $F$ cut out by $\nu_0$. Then $M$ is a biquadratic extension of $F$ containing $K$ and $M_0$. Note that
$$
{\rm Res}_K^F(\rho) \, \simeq \, {\rm Ind}_M^K(\chi) \otimes \sigma' \, \simeq \, {\rm Ind}_M^K(\chi \otimes {\rm Res}_M^K(\sigma')).
$$
Thus
$$
{\rm Ind}_M^F(\chi \otimes {\rm Res}_M^K(\sigma')) \, \simeq \, \rho \, \oplus \, \rho \otimes \delta.
$$
The left hand side is isomorphic to Ind$_{M_0}^F($Ind$_M^{M_0}(\chi \otimes {\rm Res}_M^K(\sigma')))$, which is invariant under twisting by $\nu_0$. Thus
$$
\rho \, \oplus \, \rho \otimes \delta \, \simeq \, \rho \otimes \nu_0 \, \oplus \, \rho \otimes \delta\nu_0.
$$
Since $\rho$ is irreducible, we must have
$$
\rho \, \simeq \, \rho \otimes \beta, \quad {\rm where} \quad \beta \in \{\nu_0, \nu_0\delta\}.
$$
Denote by $L$ the quadratic extension of $F$ cut out by $\beta$. (It must be either of the quadratic extensions 
which are contained in $M$ and different from $K$.) Then there exists an irreducible $2$-dimensional rpesentation $\eta$ of $\mathfrak G_L$ such that
$$
\rho \, \simeq \, {\rm Ind}_L^F(\eta),
$$
giving (ii). Lemma 3.14 is now proved.

\qed

\medskip

{\it Proof of Theorem A$^\prime$ modulo Theorem D} (contd.):

Since (i), (ii) of Lemma 3.14 coincide with (i), (ii) (respectively) of Theorem A$^\prime$, we may asume from here on that we are in case (II) (of Lemma 3.6), with $\sigma'$ non-dihedral. But as $\sigma^\theta \simeq \sigma' \otimes \mu$ for a character $\mu$, $\sigma$ is dihedral iff $\sigma'$ is. Indeed, if $\sigma$ were dihedral, it wold admit a self-twist by a quadratic character $\nu \ne 1$ and this would consequently force $\sigma'$ to admit self-twist by $\nu^\theta$. So neither $\sigma$ nor $\sigma'$ is dihedral.
We have 
$$
\sigma \, \simeq \, (\sigma^{\theta})^{\theta} \, \simeq \, (\sigma')^{\theta} \otimes \mu^{\theta}
\, \simeq \, \sigma \otimes (\mu^{\theta}/\mu).
$$
Since $\sigma$ is not dihedral, it does not admit any non-trivial self-twist by a character, 
and so we must have
$\mu = \mu^{\theta}$. So there exists a character $\nu$ of ${\mathfrak G}_F$ such that
$\mu$ is the restriction $\nu_K$ of $\nu$ to ${\mathfrak G}_K$. Thus we get (from (3.4),
$$
(\rho \otimes \nu)_K \, \simeq \, \sigma \otimes (\sigma' \otimes \mu) 
\, \simeq \, \sigma \otimes \sigma^{\theta}.
$$
It suffices to show that some character twist of $\rho$ is modular. So we may, after replacing
$\rho$ by its twist by $\nu^{-1}$, that
$$
\rho_K \, \simeq \, \sigma \otimes \sigma^{\theta}.
\leqno(3.16)
$$
If $\delta$ denotes the quadratic character of ${\mathfrak G}_F$ corresponding to $K/F$, then $\rho$
and $\rho \otimes \delta$ are the only representations for which (3.16) holds.

\medskip

It is easy to see that the induction (to ${\mathfrak G}_F$) of the exterior square of $\sigma$, 
i.e., det$(\sigma)$, is a summand of the exterior square of the induction
of $\sigma$. Thanks to semisimplicity, we may then define the 
{\it Asai representation} of $\sigma$, denoted $As(\sigma)$, by the decomposition
$$
\Lambda^2({\text Ind}_K^F(\sigma)) \, \simeq \, As(\sigma) \oplus {\text Ind}_K^F({\text det}(\sigma)).
\leqno(3.17)
$$
 
\medskip
 
\noindent{\bf Lemma 3.18} \, \it $\rho$ is isomorphic to $As(\sigma)$ or $As(\sigma) \otimes \delta$.
\rm

\medskip

{\it Proof of Lemma}. \, Let $\beta$ denote the tensor square representation of ${\text Ind}_K^F(\sigma)$,
so that
$$
\beta \, = \, \Lambda^2({\text Ind}_K^F(\sigma)) \oplus {\text sym}^2({\text Ind}_K^F(\sigma)).
$$
We can also write 
$$
\beta \, \simeq \, {\text Ind}_K^F(\sigma \otimes {\rm Res}_K^F({\text Ind}_K^F(\sigma))),
\leqno(3.19)
$$
which implies that
$$
\beta \, \simeq \, \beta \otimes \delta.
\leqno(3.20)
$$
So $As(\sigma) \otimes \delta$ must also occur in $\beta$. 

On the other hand, since the restriction to ${\mathfrak G}_K$ of 
${\text Ind}_K^F(\sigma)$ is $\sigma \oplus \sigma^{\theta}$, we get from
(3.19),
$$
\beta \, \simeq \, {\text Ind}_K^F({\text sym}^2(\sigma) \oplus \Lambda^2(\sigma)
\oplus (\sigma \otimes \sigma^{\theta})).
\leqno(3.21)
$$
Since $\rho_K$ is (by (3.16)) isomorphic to
$\sigma \otimes \sigma^{\theta}$, it must occur in the induction 
of the latter to ${\mathfrak G}_F$; ditto for the twist of
$\rho$ by $\delta$.  Hence, by the additivity of induction,
the representation on the right of (3.21) is forced to be
$$
{\text Ind}_K^F({\text sym}^2(\sigma)) \oplus 
{\text Ind}_K^F(\Lambda^2(\sigma)) \oplus \rho \oplus (\rho \otimes \delta).
$$
The lemma now follows in view of (3.17).

\qed

\medskip

So we may, after possibly replacing $\rho$ by $\rho \otimes \delta$, assume that 
$$
\rho \, \simeq \, As(\sigma),
\leqno(3.22)
$$
for an irreducible $2$-dimensional, continuous $\C$-representation 
$\sigma$ of ${\mathfrak G}_K$ with solvable image.

\medskip

For any cuspidal automorphic representation $\pi$ of GL$(2, \A_K)$, 
one may associate the following 
{\it Asai $L$-function}:
$$
L(s, \pi, r) \, = \, \prod_v \, L(s, As(\sigma_w(\pi))),
\leqno(3.23)
$$
where $v$ runs over the set $\Sigma_F$ of all the places of $F$, and
for each $v \in \Sigma_F$, $w$ denotes a place of $K$ above $v$ 
and $As(\sigma_w(\pi))$ denotes the Asai representation associated
to $\sigma_w(\pi)$. Note that the definition of $As(\sigma_w(\pi))$ is independent of the 
choice of $w$ above $v$. When $v$ splits in $K$, $K_v = K_w \times K_{\theta w}$,
and $As(\sigma_v(\eta))$
simply means the tensor product $\sigma_w(\eta) \otimes \sigma_{\theta w}(\eta)$.

The $L$-function on the left of (3.23) looks like a 
{\it Langlands 
$L$-function}, and we need to explain why we are justified in
adopting such a notation. For this recall that the 
$L$-group of the restriction of scalars of GL$(2)/K$ to $F$ is
the semidirect product
$$
{}^L(R_{K/F} {\text GL}(2)/K) \, = \, ({\text GL}(2,\C) \times 
{\text GL}(2,\C)) \times {\text Gal}(K/F),
\leqno(3.24)
$$
where $\theta$ acts by interchanging the two factors. 
One defines a representation
$$
r: {}^L(R_{K/F} {\text GL}(2)/K) \, \rightarrow {\text GL}(\C^2 \otimes \C^2) \,
\simeq \, {\text GL}(4, \C)
\leqno(3.25)
$$
by setting, for all $x, y$ in $\C^2$,
$$
r(g,g'; 1)(x \otimes y) \, = \, g(x) \otimes g(y)
$$
and
$$
r(1,1; \theta)(x\otimes y) \, = \, y \otimes x.
$$
At any finite place 
$w$ of $K$ where $\pi$ is unramified, there is a diagonal matrix
$[\alpha_w, \beta_w]$ in GL$(2, \C)$ such that
$$
L(s, \pi_w) \, = \, \frac{1}{(1-\alpha_w q_w^{-s})(1-\beta_w q_w^{-s})},
\leqno(3.26)
$$
where $q_w$ is the norm of $w$. If $v$ is any 
finite place of $F$ which is unramified for $(K/F, \pi)$, 
i.e., if $v$ is unramified in $K$ and if $\pi$ is unramified at any place $w$ 
of $K$ above $v$, then
one may associate, as in [HLR], a (Langlands) conjugacy class $A_v(\pi)$ in
${}^L(R_{K/F} {\text GL}(2)/K)$. When composed with $r$, one gets
$$
r(A_v(\pi)) \, = \, [\alpha_w \alpha_{\theta w}, \alpha_w \beta_{\theta w},
\beta_w\alpha_{\theta w}, \beta_w \alpha_{\theta w}]
\leqno(3.27)
$$
if $v$ {\it splits} into $(w, \theta w)$ in $K$, and
$$
r(A_v(\pi)) \, = \, \begin{pmatrix}
\alpha_v & 0 & 0 & 0\\
0 & 0 & \alpha_v & 0\\
0 & \beta_v & 0 & 0\\
0 & 0 & 0 & \beta_v
\end{pmatrix}
$$
if $v$ remains prime.
Since $L(s, \pi_w)$ is $L(s, \sigma_w(\pi))$, we get 
easily the identity
$$
L(s, As_v(\sigma(\pi))) \, = \, L(s, r(A_v(\pi)))
\leqno(3.28)
$$
at any finite place $v$ unramified for $(K/F, \pi)$.
This shows the appropriateness of the
notation of (3.23). It is also important because
the automorphic results we will need later will use the
Langlands formalism.

\medskip

If we admit Theorem D, we then have a
unique isobaric automorphic representation
$\Pi$ of GL$(4, \A_F)$ such that
$$
L(s, \Pi) \, = \, L(s, \pi, r).
$$
In view of the discussion above, it is clear that
$$
L(s, \Pi_v) \, = \, L(s, \rho_v)
\leqno(3.29)
$$
at almost all places $v$. By a standard argument
comparing the functional equations of $L(s, \Pi)$ and 
$L(s, \rho)$, we also get such an equality
of $L$-factors at {\it every} place $v$. 
Since $\rho_K$ is by construction associated to
$\pi \boxtimes (\pi \circ \theta)$, we get the identity
$$
L(s, \Pi_K) \, = \, L(s, \pi \boxtimes (\pi \circ \theta)).
\leqno(3.30)
$$
Since $\rho_K \simeq (\sigma \otimes (\sigma^{\theta})$
is irreducible, the main result of [Ra1] implies that
$\pi \boxtimes (\pi \circ \theta)$ is cuspidal. Since
$\Pi$ base changes to a cuspidal representation, it must itself be 
cuspidal by [AC]. 

This finishes the proof of Theorem A$^\prime$, and hence Theorem A, 
modulo Theorem D.

\qed

\vskip 0.2in

\section
{\bf Reduction to weak lifting, and the cuspidality criterion}

\bigskip

We will begin the proof of Theorem D in this section and finish in section 6. 

Fix $K/F$ quadratic with non-trivial automorphism $\theta$ as above, and an arbitrary 
cuspidal
automorphic representation $\pi$ of GL$(2, \A_K)$. We have to show that there exists an 
isobaric automorphic representation $\Pi$ of GL$(4, \A_F)$ such that 
$$
L(s, \Pi) \, = \, L(s, \pi; r).
$$
Suppose $\chi$ is an idele class character of $K$ with restriction $\chi_0$ to $F$. Recall 
that $\chi_0$ corresonds to the {\it transfer} from $K$ to $F$ of the character of the Weil 
group $W_K$ associated to $\chi$ by class field theory. (By abuse of notation we will use 
the same letter to signify the characters of $\A_K^\ast/K^\ast$ and $W_K$.) At any place $v$ 
of $F$ with divisor $w$ in $K$, we then get, from the definition of the Asai representation,
$$
As(\sigma(\pi_w \otimes \chi_w)) \, \simeq \, As(\sigma(\pi_w)) \otimes \chi_{0,v}.
$$
Consequently,  
$$
L(s, \pi \otimes \chi; r) \, = \,  L(s, \pi; r \otimes \chi_0).
\leqno(4.0)
$$
One knows, cf. [HLR], that $L(s, \pi; r)$ admits a meromorphic continuation to the whole 
$s$-plane and satisfies a standard functional equation.

\medskip

\noindent{\bf Proposition 4.1} \quad \it Let $\pi$ be
a cuspidal automorphic representation of GL$(2, \A_K)$. Suppose we have constructed a weak Asai 
lifting, i.e.,
an isobaric 
automorphic representation $\Pi$ of GL$(4, \A_F)$
satisfying the following identity at almost all places $v$ of $F$:
$$
L(s, \Pi_v) \, = \, L(s, \pi_w; r),
\leqno(L_v)
$$
where $w$ is any place of $K$ above $v$.
Then we have ($L_v$) and ($\varepsilon_v$) at {\bf all} the finite 
places $v$, and ($L_\infty$) as well.
In addition, the cuspidality criterion of Theorem D holds. 
\rm

\medskip

{\it Proof}. \, 
Let $S$ be the (finite) set of places of $F$ outside which ($L_v$) holds.
Note first that the central character $\Omega$ of $\Pi$ is simply the restriction $\omega_0$ to $F$ 
of the central character $\omega$ of $\pi$. Indeed at any $v$, $\omega_0$ corresponds to the transfer 
of the Galois character attached to $\omega$, which also gives the determinant of $As(\sigma_w(\pi))$, 
with $w$ being a place of $K$ above $v$. Since by definition, $L(s, \pi_w; r)$ equals 
$L(s, As(\sigma_w(\pi))$, it follows that the idele class
characters $\Omega$ and $\omega_0$ agree outside $S$, and hence agree everywhere by a classical 
result of Hecke.

\medskip

Now some notations. If $f(s), g(s)$ are two meromorphic functions of
$s$ such that their quotient is invertible, we will write
$f(s) \, \equiv \, g(s)$. At any place $v$, given a character $\nu$ of 
$F_v^\ast$, we can write it as $\nu_0|.|^z$, for a unitary character
$\nu_0$ and a complex number $z$. The real part of $z$ is uniquely 
defined; we will
call it the {\it exponent of $\nu$}, and denote it $e(\nu)$.

\medskip

Choose a finite order character $\mu$ of $C_F$ with $\mu_\infty = 1$ 
(which means that $\mu$ is totally even)
such that $\mu_v$ is sufficiently ramified at every finite place $u$
in $S$ so as to make the $L$-factors at $u$ of $\Pi \mu_u$, $(\pi, r \otimes \mu_u)$,
and their contragredients, all equal $1$. This is 
evident for $L(s, \pi; r \otimes \mu)$ as its local factors are defined here to coincide with the 
corresponding Galois factors, and it is possible for $L(s, \Pi \otimes \mu)$ by the results of 
[JPSS1]. Comparing the global functional equations of both $L$-functions,  and noting that twisting by
$\mu$ does not change anything at infinity, 
we get 
$$
\prod\limits_{v | \infty} \, L(s, \Pi_v) 
L(1-s, \pi_w^\vee; r)) \, \equiv
\, \prod\limits_{v | \infty} \, L(1-s, \Pi_v^\vee) 
L(s, \pi_w; r).
\leqno(4.2)
$$
For {\it any} place $v$, archimedean or otherwise, for any $n \geq 1$,
and for any
cuspidal automorphic representation 
$\eta$ of GL$(n, \A_F)$, it is known that
$L(s, \eta_v)$ is holomorphic in $\Re(s) > \frac{1}{2} - t$, for some
$t = t(\eta, v) > 0$ (see [BaR], Prop. 2.1, part B). Consequently, 
$L(s, \Pi_v)$ has no pole in common with 
$L(1-s, \Pi_v^\vee)$. 

Thus the poles of $L(s, \Pi_\infty)$ are contained in those of $L(s, \pi_\infty; r)$. For the 
converse direction, we appeal to the factorization formula
$$
L(s, \pi \times (\pi \circ \theta)) \, = \, L(s, \pi; r)L(s, \pi; r \otimes \delta),
\leqno(4.3)
$$
where $\delta$ is the quadratic character of $F$ associated to $K/F$. This formula is evident 
from the definition (3.23). Since the local factors are never zero at any place $v$, the 
$v$-factor of $L(s, \pi; r)$ can have a pole somewhere only if the $v$-factor of 
$L(s, \pi \times (\pi \circ \theta))$ also does. On the other hand, by Theorem M of [Ra1], 
there is a unique isobaric automorphic representation $\pi \boxtimes (\pi \circ \theta)$ 
of GL$(4, \A_F)$ whose standard $L$-function coincides with $L(s, \pi \times \pi')$. So, 
applying [BaR], Proposition 2.1, again, we see that $L(s, \pi_v; r)$ is holomorphic in 
$\Re(s) > \frac{1}{2} - t$, for some $t > 0$. Hence the poles of $L(s, \pi_\infty; r)$ 
are distinct from those of $L(1-s; \pi^\vee_\infty; r)$, and so must coincide with those
of $L(s, \Pi_\infty)$. Let us use the customary notation $\Gamma_\R(s) = \pi^{-s/2}\Gamma(\frac{s}{2})$ and 
$\Gamma_\C(s) = 2(2\pi)^{-s}\Gamma(s)$. Using the duplication formula expressing
$\Gamma_\C(s)$ as $\Gamma_\R(s)\Gamma_\R(s+1)$, and appealing to the fact that $L(s, \Pi_\infty)$
is a standard $L$-factor of GL$(4, F_\infty)$ and that $L(s, \pi_\infty; r)$ is (by definition) a Galois $L$-factor,
we may write
$$
L(s, \Pi_\infty) \, = \, \prod_{j=1}^n \Gamma_\R(\frac{s+a_j}{2}) \quad {\rm and} \quad
L(s, \pi_\infty; r) \, = \, \prod_{j=1}^n \Gamma_\R(\frac{s+b_j}{2}),
$$  
for some complex numbers $a_j, b_j, 1 \leq j \leq n = 4[F :\Q]$. (The fact that $n$ is $4$ times $[F:\Q]$ will play no role.) We may renumber the $a_j$ and $b_j$ and assume that there exists an integer $m$, with $0 \leq m \leq n$, such that $a_j = b_j$ for all $j > m$ and the sets $\{a_j \vert j \leq m\}$ and
$\{b_j \vert j \leq m\}$ are totally disjoint. We have nothing to prove if $m = 0$, so assume that $m$ is positive. Now we appeal to the following 

\medskip

\noindent{\bf Baby Lemma} \, Let $m > 0$ be an integer and let $\{a_j \vert j \leq m\}$, 
$\{b_j \vert j \leq m\}$ be subsets of $\C$ with empty intersection. Then the polar divisors of 
$L_1(s): = \prod_{j=1}^n \Gamma_\R(\frac{s+a_j}{2})$ and
$L_2(s): = \prod_{j=1}^n \Gamma_\R(\frac{s+b_j}{2})$ cannot be the same.
\rm

\medskip

{\it Proof of Baby Lemma}. \, Clearly the support of the polar divisor of 
$\Gamma_\R(\frac{s+c}{2})$ is, for any $c \in \C$, the set 
$\{-c-2k \vert k \in \Z, k \geq 0\}$. 
Suppose the 
poles of $L_1(s)$ and $L_2(s)$ coincide. Then $-a_1$ must be a pole of some 
$\Gamma_\R(\frac{s+b_j}{2})$. After renumbering the $b_j$, 
we may then assume that $a_1 = b_1 + 2k_1$ 
for some positive integer $k_1$. (Since $a_1 \ne b_1$, $k_1$ cannot be $0$.) 
Now $-b_1$ will need to be a pole of some $\Gamma_\R(\frac{s+a_j}{2})$, and
$j$ must be
$\geq 2$. After renumbering the $a_j$ 
we may assume that $b_1 = a_2 + 2\ell_1$ for some $\ell_1 > 0$. We may continue thus and write, 
after suitable renumberings of the $a_j$ and the $b_j$, $a_2 = b_2 + 2k_2$, $b_2 = a_3 + 2\ell_2$,
and so on. This leads to the string of inequalities $a_1 < b_1 < a_2 < b_2 < a_3 \dots < a_m < b_m$. 
Then $-b_m$ is a pole of $L_2(s)$ and it is not a pole of $L_1(s)$.
\qed

The identity $(L_\infty)$ now folows.

\medskip

Next fix a finite place $v$ in $S$. Now choose a character $\mu$ which is unramified at $v$, 
but is highly ramified at every $u$ in $S - \{v\}$. Comparing the functional equations and 
arguing excatly as in the archimedean case, we get the desired identity $(L_v)$.  

\medskip

Next we prove the identity
of epsilon factors. Fix any $v$ in $S$
and note that by 
[JPSS2], for $\mu_u$ sufficiently ramified at $u \in S - \{v\}$, the epsilon factor of
$\Pi_u \otimes \mu_u$ depends only on 
$\mu_u$ and 
$\omega_{0,u}$, and the dependence is simple. 
Similarly, by [DeH], the epsilon factor at $u$ of
$(\pi,r \otimes \mu)$ has the same dependence on $\mu_u$
and $\omega_{0,u}$; the reason we can apply [DeH] is that we have  defined
the local factors of $(\pi,r \otimes \mu)$ as those associated to the corresponding
representations of the local Weil groups. The analogous 
statements hold for the contragredients, and this results in
the identity  
($\varepsilon_{v}$) as we already know that the $L$-factors agree. 
This finishes the proof of the first part of Proposition 4.1.

\bigskip

It is left to prove the {\bf cuspidality criterion}. 

Note that (4.3) implies the identity
$$
L(s, \pi \times (\pi \circ \theta)) \, = \, L(s, \Pi_K),
\leqno(4.4)
$$
where $\Pi_K$ denotes the base change of the isobaric representation $\Pi$ to $K$. 
This forces, by the existence and uniqueness of $\pi \boxtimes (\pi \circ \theta)$ 
(cf. Theorem M of [Ra1]), we see that
$$
\Pi_K \, \simeq \, \pi \boxtimes (\pi \circ \theta).
\leqno(4.5)
$$

First suppose $\Pi$ is cuspidal. If $\Pi_K$ is not cuspidal, then by the theory of base 
change ([AC]), $\Pi$ must be automorphically induced from a cuspidal automorphic 
representation $\eta$ of GL$(2, \A_K)$; write $\Pi = I_K^F(\eta)$. Then we must have
$$
\Pi_K \, = \, \eta \boxplus (\eta \circ \theta).
\leqno(4.6)
$$
The idea now is to compute the exterior square $L$-function
of $\Pi_K$ in two different ways. We refer to [JS] and [BF] for the relevant facts about 
these degree $6$ $L$-functions. On the one hand, 
(4.6) gives
$$
L^S(s, \Pi_K, \Lambda^2) \, = \, 
L^S(s, \eta \times (\eta \circ \theta))L^S(s, \omega_\eta)L^S(s, \omega_\eta^\theta),
\leqno(4.7)
$$
where $\omega$ (resp. $\omega_\eta^\theta$) is the central character of $\eta$ 
(resp. $\eta \circ \theta$). This can be
seen easily at the unramified places $w$. Indeed, the above identity is induced 
by the following:
$$
\Lambda^2(\sigma_w(\eta) \otimes \sigma_w(\eta \circ \theta)) \, \simeq \, 
\sigma_w(\eta) \otimes \sigma_w(\eta \circ \theta) \, \oplus \, {\rm det}(\sigma_w(\eta)) \,
\oplus \, {\rm det}(\sigma_w(\eta \circ \theta)),
$$
which is easy to verify.
Consequently, $L^S(s, \Pi_K, \Lambda^2)$ is divisible
by two abelian $L$-functions, namely $L^S(s, \omega_\eta)$ and 
$L^S(s, \omega_\eta^\theta)$. 
On the other hand, we also have the identity
$$
\Lambda^2(\sigma_w(\pi) \otimes \sigma_w(\pi \circ \theta)) \, \simeq \, 
{\rm det}(\sigma_w(\pi)) \otimes {\rm sym}^2(\sigma_w(\pi \circ \theta) \, \oplus \, 
{\rm sym}^2(\sigma_w(\pi)) \otimes {\rm det}(\sigma_w(\pi \circ \theta)),
$$
implying the following equality of $L$-functions:
$$
L^S(s, \Pi_K, \Lambda^2) \, = \, 
L^S(s, {\rm sym}^2(\pi) \otimes \omega)
L^S(s, {\rm sym}^2(\pi') \otimes (\omega \circ \theta)).
\leqno(4.8)
$$

Suppose $\pi$ is {\bf non-dihedral}. Then it is known (cf. [GeJ]) that sym$^2(\pi)$ is
cuspidal; so is sym$^2(\pi \circ \theta)$. Consequently, thanks to
(4.8), $L^S(s, \Pi_K, \Lambda^2)$
cannot be divisible by an abelian $L$-function, leading to a contradiction, and so 
$\Pi_K = \pi \boxtimes (\pi \circ \theta)$ must be cuspidal for non-dihedral $\pi$. 
In this case, if $\pi \circ \theta$ were isomorphic to $\pi \otimes \chi$ for a 
character $\chi$, then $L(s, \Pi_K \otimes (\chi\omega)^{-1})$ would have a pole at 
$s=1$, contradicting the cuspidality of $\Pi_K$. So the cuspidality of $\Pi$ implies, 
when $\pi$ is non-dihedral, that $\pi \circ \theta$ cannot be a character twist of 
$\pi$, as asserted in Theorem D.

Conversely, suppose that $\pi$  is non-diheedral and 
{\it not} equivalent to any character twist of $\pi \circ \theta$. Then by the 
cuspidality criterion in Theorem M of [Ra1], we know that $\pi \boxtimes (\pi \circ 
\theta)$ is a cuspidal automorphic representation of GL$(4, \A_K)$.  But this is 
just $\Pi_K$ by (4.5). Hence $\Pi$ base changes to a something cuspidal over K, and hence must be cusspidal  itself (cf. [AC]). So the cuspidality criterion of  Proposition 4.1
is now  proven in the {\it non-dihedral} case. 

\medskip

Now let $\pi$ be {\bf dihedral}, so that 
$\pi$ is an automorphically induced representation $I_M^K(\chi)$ 
for an idele class character $\chi$ of a quadratic extension 
$M$ of $K$. Denote the corresponding $2$-dimenssional
representation of the Weil group $W_K$ by $\sigma = $Ind$_M^K(\chi)$. 
Since $\pi$ is cuspidal, $\sigma$ s irreducible. In this (dihedral) case
there is a $4$-dimensional global Asai representation 
$As(\sigma)$, of $W_F$ such that for any character $\nu$  of $F$, we have
$$
L(s, \pi; r) \, = \, L(s, As(\sigma)),
\leqno(4.9)
$$
and by the Tchebotarev density theorem,
$$
As(\sigma)\vert_{W_K} \, \simeq \, \sigma \otimes 
\sigma^{\theta}.
\leqno(4.10)
$$

\medskip

\noindent{\bf Lemma 4.11} \, \it Let $\pi$ be dihedral with associated 
representation $\sigma = $Ind$_M^K(\mu)$ of $W_K$, where M/K is
quadratic and $\mu$ a character of $W_M$. Then 
$\Pi$ is cuspidal iff $As(\sigma)$ is 
irreducible.
\rm

\medskip

{\it Proof}.  Suppose  $As(\sigma)$ is reducible. It is known that the 
$\C$-representations of $W_F$ are completely reducible. So we  may write $As(\sigma)
= \oplus_j n_j \eta_j$ with eachh $n_j \geq 1$ and
$\eta_j$ a proper irreducible summand of  $As(\sigma)$, with  $\eta_i, \eta_j$ 
inequivalent if $i \ne j$. Then it is elementary to see that
$$
{\rm dim}_\C {\rm Hom}_{W_F}(1, As(\sigma)  \otimes As(\sigma)^\vee) 
\, = \, \sum_j n_j^2.
\leqno(4.12)
$$
On the other hand, since $As(\sigma)$ is a $\C$-representation of $W_F$, 
we may use Brauer's theorem ([De]) and get a virtual sum decomposition
$$
As(\sigma)  \otimes As(\sigma)^\vee \, \simeq \, \oplus_{i=1}^r m_i
{\rm Ind}_{L_i}^F(\lambda_i),
$$
where $m_i$ is, for each $i \leq r$, an integer, $L_i$ a finite extension of $F$,
and $\lambda_i$ a character of $W_{L_i}$. By the inductivity and additivity of
$L$-functions, we see that
$$
L(s, As(\sigma) \otimes As(\sigma)^\vee) \, = \, \prod_{i=1}^r L(s, \lambda_i)^{m_i}.
$$
Moreover one knows by Hecke that $L(s, \lambda_i)$ is invertible at $s=1$ unless
$\lambda_i$ is the trivial character, in which case it has a pole
of order $1$. From this one gets
$$
-{\rm ord}_{s=1} L(s, As(\sigma) \otimes As(\sigma)^\vee) \, = \, 
{\rm dim}_\C {\rm Hom}_{W_F}(1, As(\sigma)  \otimes As(\sigma)^\vee)
\leqno(4.13)
$$
Combining (4.12) and (4.13) we see that 
$L(s, \Pi \times  \Pi^\vee)$ has a pole of order $> 2$ at $s=1$, which can only 
happen if $\Pi$ is non-cuspidal ([JS]). 
The converse assertion holds by reversing the argument.

\qed

\medskip

{\it Proof of the cuspidality criterion} (contd.):

Fix an extension of $\theta$ to the Galois closure $\tilde M$
of $M$ over $F$, and denote it again by $\theta$. Denote by 
$\alpha$ the non-trivial automorphism of $M/K$, and by $\epsilon$ the
quadratic character of $W_K$ corresponding to $M/K$. Then $\epsilon^\theta$
is the quadratic character of $W_F$ corresponding to $M^\theta \ne M$.
Note that
$$
\sigma^\theta \, \simeq \, {\rm Ind}_{M^\theta}^K(\chi^\theta),
\leqno(4.14)
$$
where $\chi^\theta$ is the character of $W_{M^\theta}$ defined by
transporting $\chi$ via $\theta$.

We have to show now that $As(\sigma)$ 
is irreducible iff $M/F$ is non-Galois {\it and} the representation
$$
\tau: = \, {\rm Res}_M^K(\sigma^{\theta}) \otimes \chi
\leqno(4.15)
$$
of $W_M$ does not extend to $W_F$. 

Suppose $M/F$ is Galois. We have by Mackey and the 
definition of $As(\sigma)$,
$$
{\rm Res}_K^F(As(\sigma)) \, \simeq \, \sigma \otimes \sigma^{\theta} 
\, \simeq \, {\rm Ind}_M^K(\chi\chi^\theta) \oplus 
{\rm Ind}_M^K(\chi\chi^{\theta\alpha}).
\leqno(4.16)
$$
The first summand on the right is evidently $\theta$-invariant, and so
extends to a $2$-dimensional representation of $W_F$ occurring in $As(\sigma)$.
Hence $As(\sigma)$ is reducible in this case.

\medskip

So we may assume from here on that $M/F$ is non-Galois.
We have by the definition of $\tau$,
$$
\sigma \otimes \sigma^{\theta} \, \simeq \,
{\rm Ind}_M^K(\tau).
\leqno(4.17)
$$
Restricting to $W_M$ we get
$$
{\rm Res}_M^F(As(\sigma)) \, \simeq \, {\rm Res}_M^K(\sigma \otimes
\sigma^{\theta}) \, \simeq \, {\rm Res}_M^K({\rm Ind}_M^K(\tau)) \, \simeq \,
\tau \oplus \tau^{\alpha}.
\leqno(4.18)
$$

\medskip

First suppose that $\tau$ extends to $F$. Already the fact that it extends 
to a representation $\tau^K$, say, of $W_K$ implies that 
$\sigma \otimes \sigma^{\theta}$ is reducible. 
More precisely,
$$
\sigma \otimes \sigma^{\theta} \, \simeq \, \tau^K \oplus
(\tau^K \otimes \epsilon).
\leqno(4.19)
$$
Now the existence of an extension of $\tau$ all the way to $F$ says  
that $\tau^K$ extends to a $W_F$-summand of dimension $2$ 
in $As(\sigma)$. Hence $As(\sigma)$ is reducible.

\medskip

Now we will assume that $As(\sigma)$ is reducible, and prove the 
more subtle {\it converse}. Then $\sigma \otimes \sigma^{\theta}$ is reducible,
and by (4.17), either $\tau$ is reducible or $\tau$ extends to $K$.

First consider when {\it $\tau$ is irreducible}, but extends to a 
representation $\tau^K$ of $W_K$. Then (4.19) will hold, and moreover, 
the $\theta$-invariance of $\sigma \otimes \sigma^{\theta}$ implies 
one of the following:

\noindent(4.20)
\begin{itemize}
\item[{(i)}] $\tau^K$ is $\theta$-invariant, i.e., $\tau$ extends to $F$;
\item[{(ii)}] $(\tau^K)^{\theta} \, \simeq \, \tau^K \otimes \epsilon$.
\end{itemize}
Suppose (ii) happens without (i). Then
$$
(\tau^K)^{\theta} \, \simeq \, \tau^K \otimes \epsilon
$$
and this implies by combining (4.17) and (4.19), that $As(\sigma)$ 
is irreducible, a contradiction! So (i) must hold and 
we are done if $\tau$ is irreducible.

\medskip

It is left to consider when $\tau$ is {\it reducible}. 
For this case we need 
the following

\medskip

\noindent{\bf Lemma 4.21} \, \it Let $\tau$ be reducible, with $M/F$ non-Galois. 
Then there exists 
a character $\mu$ of $W_M$ such that
\begin{enumerate}
\item[{(a)}] \, $\sigma^{\theta} \, \simeq \, {\rm Ind}_M^K(\mu)$;
\item[{(b)}] \, $\tau \, \simeq \, (\mu \oplus \mu^\alpha) \otimes \chi$;
\item[{(c)}] \, ${\rm Res}_K^F(As(\sigma)) \, \simeq \, {\rm Ind}_M^K(\mu\chi) 
\oplus {\rm Ind}_M^K(\mu^\alpha \chi)$; \, and
\item[{(d)}] \, We have
$$
\chi/\chi^\alpha \, = \, \mu/\mu^\alpha \, = \, \epsilon^{\theta}_M.
$$
\end{enumerate}
\rm

\medskip

{\it Proof of Lemma}. 

(a) \, $\tau$ is reducible iff Res$_M^K(\sigma^{\theta})$ is 
reducible, which happens, thanks to the irreducibility of $\sigma^{\theta}$, 
iff $\sigma^{\theta}$ is Ind$_M^K(\mu)$ for some character $\mu$ of $W_M$.

(b) By Mackey, Res$_M^K({\rm Ind}_M^K(\mu)) \, \simeq \, \mu \oplus \mu^\alpha$.
So by the definition (4.14), $\tau \simeq (\mu \oplus \mu^\alpha) \otimes \chi$.

(c) In view of (4.17), this part follows from (b).

(d) Since $M/F$ is non-Galois, the restriction $\epsilon^\theta_M$ of $\epsilon^\theta$
to $W_M$ is non-trivial. And as
$\sigma^{\theta}$ is induced from $M^\theta$, we have $\sigma^{\theta}$ 
is isomorphic to $\sigma^{\theta} \otimes \epsilon^{\theta}$. Then (4.14) gives
the isomorphism $\tau \, \simeq \, \tau \otimes \epsilon^{\theta}_M$. Using (b), we get
$\mu\chi \oplus \mu^\alpha\chi \, \simeq \, \mu\chi\epsilon^\theta_M 
\oplus \mu^\alpha\chi\epsilon^\theta_M$. This forces the desired equality 
$\mu^\alpha/\mu = \epsilon^\theta_M$, which shows in particular that 
$\mu^\alpha/\mu$ extends to $K$. On the other hand, since $\sigma^\theta$ is 
induced from $M$, it follows that $\sigma$ is induced from $M^\theta$, and so
we must have $\sigma \otimes \epsilon^\theta \, \simeq \, \sigma$. In other words,
Ind$_M^K(\chi\epsilon^\theta_M) \, \simeq \, {\rm Ind}_M^K(\chi)$. Restricting to $M$,
we get the isomorphism $\chi\epsilon^\theta_M \oplus \chi^\alpha\epsilon^\theta_M
\, \simeq \, \chi \oplus \chi^\alpha$, which gives $\chi^\alpha/\chi 
= \epsilon^\theta_M$. Done.

\qed

\medskip

{\it End of proof of the cuspidality criterion}:

Using part (d) of the
Lemma above, we see that 
$$
(\mu\chi)^\alpha \, = \, (\mu^\alpha/\mu)(\chi^\alpha/\chi)\mu\chi \, 
= \, (\epsilon^\theta_M)^2\mu\chi \, = \, \mu\chi.
$$
Hence $\mu\chi$ extends to a character $\nu$, unique up to multiplication
by $\epsilon$, of $W_K$, and consequently,
Ind$_M^K(\mu\chi) = \nu \oplus \nu\epsilon$.
Then $\mu^\alpha\chi = \epsilon^\theta_M\mu\chi$ can be extended to $K$, either
as $\epsilon^\theta\nu$ or as $\epsilon\epsilon^\theta\nu$.
In any case, part (c) of
the Lemma gives
$$
{\rm Res}_K^F(As(\sigma)) \, \simeq \, \nu \oplus \nu\epsilon \oplus
\nu\epsilon^\theta \oplus \nu\epsilon\epsilon^\theta.
$$
Since Res$_K^F(As(\sigma))$ is $\theta$-invariant, we need
$$
\nu^\theta \, \in \, \{\nu, \nu\epsilon, \nu\epsilon^\theta, 
\nu\epsilon\epsilon^\theta\}.
$$
Clearly $\nu^\theta$ cannot be $\nu\epsilon$ or $\nu\epsilon^\theta$
as it will contradict the identity $(\nu^\theta)^\theta = \nu$. So
$\nu^\theta$ must be $\nu$ or $\nu\epsilon\epsilon^\theta$. In either case, 
take the extension $\nu\epsilon\epsilon^\theta$ of $\mu^\alpha\chi$,
so that the representation 
$$
\tau^K: = \, \nu \oplus \nu\epsilon\epsilon^\theta
$$
of $W_K$ extends $\tau$ and more importantly, extends to a representation of
$W_F$, as we needed to show.

\qed

\vskip 0.2in

\section
{\bf Distinguished representations} 
 
\bigskip

Let $K/F$ be a quadratic extension of number fields with
Gal$(K/F) = \{1, \theta\}$. The object of this section is to 
establish Theorem D for the nice subclass of {\it distinguished}
cusp forms $\pi$ ([HLR]) on GL$(2)$/K. It is necessary to treat this case separately as 
certain twists of the Asai $L$-function of $\pi$ will, in such a case, 
admit poles, complicating the argument using the converse theorem, 
which we will utilize for $\pi$ of {\it general type} in
the next section.

\medskip   

We will use the following notation. If $\chi$ is an idele class character
of $K$, we will write $\chi_0$ for its
restriction to $F$. (This corrsponds to taking the {\it transfer} of the 
associated Galois character.) 
Moreover, if $\mu$ is a character of $F$, then we will write $\mu'$ to 
signify any character of $K$ such that $\mu = \mu'_0$. If $\mu_1'$ is
another extension of $\mu$, then there exists a character $\nu$ of $K$ such that
$$
\mu_1' \, = \, \mu'(\nu/(\nu \circ \theta)).
\leqno(5.1)
$$
This is because any character of $K$ whose restriction to $F$
is trivial lies in Ker$(\theta -1)$.

\medskip

Let $\pi$ be a
cuspidal automorphic representation of GL$(2, \A_K)$
with space ${\mathcal V}_\pi$. If $\mu$ is a unitary character of $F$,
then $\pi$ is said to be {\it $\mu$-distinguished} ([HLR]) 
iff the following {\it $\mu$-period integral} is non-zero for some 
function $f$ in ${\mathcal V}_\pi$: 
$$
{\mathcal P}_\mu(f) :\, = \, \int_{H(F)Z_H(F_\infty)^+\backslash H(\A_F)} 
\mu({\text det}(h))f(h) dh,
\leqno(5.2)
$$
where $H$ denotes GL$(2)/F$ with center $Z_H$, and $dh$ is the 
quotient measure induced by the Haar measure on $H(\A_F)$.
It may be useful to note for the uninitiated that when $F = \Q, K$ real quadratic,
and $f \in \pi$ a holomorphic newform of weight $(2,2)$, ${\mathcal P}_\mu(f)$ is
the $\mu$-twisted integral of the $(1,1)$ differential form $(2\pi i)^2f(z_1, z_2)
dz_1 \wedge \overline{dz}_2$ on the associated Hilbert modular surface
over (the homology class of) the modular curve; so one is justified in calling this
a period integral. 

A basic result of [HLR], section 2, asserts that, once we have fixed an
extension $\mu'$ of $\mu$, the necessary and sufficient condition for $\pi$ to be
$\mu$-distinguished is that there exists a cuspidal automorphic 
representation $\pi_0$ of $H(\A_F)$ with central character $\nu\delta$ such that
$$
\pi_{0,K}  \, \simeq \, \pi \otimes \nu'\mu',
\leqno(5.3)
$$
for a suitable extension $\nu'$ of $\nu$. 
 
\medskip

Fix such a $\mu$-distinguished $\pi$ with $(\pi_0, \nu)$ as above.
Since $\pi \boxtimes (\pi \circ \theta)$ is $\theta$-invariant,
it descends (by [AC]) to an isobaric automorphic representation 
of GL$(4, \A_F)$. We can give an explicit candidate for this descent by setting
$$
\Pi \,: = \, {\text sym}^2(\pi_0) \otimes \delta(\mu\nu)^{-1} \,
\boxplus \, \delta\mu^{-1}.
\leqno(5.4)
$$
That the base change $\Pi_K$ is $\pi \boxtimes (\pi \circ \theta)$ is easily
deduced from (5.3). There are at least four possible descents, namely
by leaving in or removing the character $\delta$ at the places where it appears in (5.4),
and this is why we needed to make a specific choice.
Note also that the automorphic induction of $\pi$ to $F$ satisfies
$$
I_K^F(\pi) \, \simeq \, \pi_0 \boxtimes I_K^F((\mu'\nu')^{-1}).
\leqno(5.5)
$$

\medskip

It suffices, by Proposition 4.1, to prove that the local factors of 
$L(s, \Pi)$ and $L(s, \pi; r)$ agree
almost everywhere.
Let $v$ be a finite place where $\pi$ and $K/F$ are unramified.
If $v$ splits in $K$, the desired identity is immediate. So assume $v$ is inert, and denote the
unique place of $K$ above it by $w$. Recall that
the exterior square of a tensor product $V \otimes W$ is
the direct sum of
sym$^2(V) \otimes \Lambda^2(W)$ and 
sym$^2(W) \otimes \Lambda^2(V)$. 
Using this conjunction with (5.5), and by the
compatibility of local and global automorphic
induction, we have
$$
\Lambda^2(\sigma_v(I_K^F((\pi))) \, \simeq \, 
{\text sym}^2(\sigma_v(\pi_0)) \otimes \delta_v(\mu_v\nu_v)^{-1} \,
\oplus {\text sym}^2({\text Ind}_{K_w}^{F_v}((\mu'_w\nu'_w)^{-1})) \otimes \nu_v\delta_v.
\leqno(5.6)
$$  
We also have
$$
{\text sym}^2({\text Ind}_{K_w}^{F_v}((\mu'_w\nu'_w)^{-1})) \, \simeq \,
{\text Ind}_{K_w}^{F_v}((\mu'_w\nu'_w)^{-2}) \oplus (\mu_v\nu_v)^{-1}.
\leqno(5.7)
$$
Combining these two identities with the fact that the induced module on the right
of (5.7) is simply the induction of the determinant of $\sigma_w(\pi)$, we get,
from the definition of the Asai representation
$$
As_v(\sigma(\pi)) \, \simeq \, {\text sym}^2(\sigma_v(\pi_0)) \otimes \delta_v(\mu_v\nu_v)^{-1} \,
\oplus \delta_v\mu_v^{-1}.
\leqno(5.8)
$$
Its $L$-factor, in view of (5.4), coincides with that of $\sigma_v(\Pi)$. Done.

\qed

\vskip 0.2in

\section
{\bf Twisted Asai $L$-functions}

\bigskip

Let $K/F$ be a quadratic extension of number fields, and let $\pi$
be a cuspidal automorphic representation of GL$(2, \A_K)$ of 
central character $\omega = \omega_\pi$.

\medskip

Let $m \in \{1,2\}$. Then for any
cuspidal automorphic representation $\eta$ of GL$(m, \A_F)$, one may define 
the {\it $\eta$-twisted Asai $L$-function} of $\pi$ by setting
$$
L(s, \pi; r \otimes \eta) \, = \, \prod_v L(s, As(\sigma(\pi_w)) \otimes \sigma_v(\eta)),
\leqno(6.1)
$$
where for each $v$ we have chosen a place $w$ of $K$ above it. 
This $L$-function converges normally in a right half plane and defines an invertible
holomorphic function there. Of course, when $m=1$, $\eta$ is simply
an idele class character of $F$ with contragredient $\eta^\vee = \eta^{-1}$.

For any idele class character $\chi$ of $K$ with restriction $\chi_0$ to $F$, we get
$$
L(s, \pi \otimes \chi; r \otimes \eta) \, = \,
L(s, \pi; r \otimes (\eta \otimes \chi_0)).
\leqno(6.2)
$$
It also follows from the definition that
$$
L(s, \pi \boxtimes (\pi \circ \theta) \otimes \eta_K) \, = \,
L(s, \pi; r \otimes \eta)L(s, \pi; r \otimes (\eta \otimes \delta)),
\leqno(6.3)
$$
where again $\delta$ denotes the quadratic character of $F$ defined by $K/F$. 

\medskip

The object of this section 
is to establish, under some loal hypotheses, the
needed analytic properties of these $\eta$-twisted Asai $L$-functions. We need the
following

\medskip

\noindent{\bf Proposition 6.4} \, \it Let $F$ be a totally imaginary number field,
$K/F$ a quadratic extension with associated character $\delta$ of $W_F$ and 
non-trivial automorphism $\theta$ of $K/F$. Let $S$ be a non-empty
finite set of
finite places of $F$ which split
in $K$, and $\pi$ a non-distinguished, cuspidal
automorphic representation of GL$(2, \A_K)$, 
which is unramified at any finite place not above $S$. Assume 
moreover that the square of the central character
$\omega_\pi$ is ramified at some place in $S$. 
Let $\eta$ be a cuspidal automorphic representation of
GL$(m, \A_F)$, $m = 1,2$, which is 
unramified at any finite place in $S$. Then we have
the following:
\begin{enumerate}
\item[(MC)] \, $L(s, \pi; r \otimes \eta)$ admits a meromorphic 
continuation
to the whole $s$-plane;
\item[(FE)] \, There is a functional equation
$$
L(1-s, \pi^\vee; r \otimes \eta^\vee) \, = \, 
\varepsilon(s, \pi; r \otimes \eta)
L(s, \pi; r \otimes \eta);
$$
\item[(E)] \, $L(s, \pi; r \otimes \eta)$ is entire; \quad and
\item[(BV)] \, $L(s, \pi; r \otimes \eta)$ is bounded in 
vertical strips of 
finite width.
\end{enumerate}
\rm

\medskip

\noindent{\bf Remark 6.5}: \, 

(a) \, When we say that $\pi$ is non-distinguished, we mean that it is not 
$\mu$-distinguished for
{\it any} character $\mu$ of $F$; so being distinguished (or not)
is a property shared
by all the character twists. 

(b) \, We will get this without any hypotheses in the next section {\it after} 
completing the proof of Theorem D.

\medskip

\noindent{\it Proof}. \. Let $\pi, \eta$ be as in the Proposition.
Let $\phi = \phi_1 \otimes \phi_2$, with $\phi_1$, resp. $\phi_2$,
lying in the space of $\pi$, resp. $\eta$. Denote by $\omega$ the
restriction of $\omega_\pi$ to $F$ times the central character
$\omega_\eta$ of $\eta$.

{\it (FE) and (MC)}: \,
There is a closely related $L$-function of 
the pair $(\pi, \eta)$, which we will denote by
$L_1(s, \pi; r \otimes \eta)$, given by an {\it integral
representation}. For $m=1$, this was done in 
the work of Harder, Langlands and Rapoport ([HLR]), generalizing
the earlier construction of Asai for holomorphic Hilbert modular forms.
For $m=2$, which is the more difficult case, 
this was done in the work of Piatetski-Shapiro and Rallis
([PS-R]). (Their work was motivated by the earlier work of 
P. Garrett on the triple product $L$-function attached to triples
of cusp forms on GL$(2)$, and there is a formal similarity between such
$L$-functions and GL$(2)$-twisted Asai $L$-functions.)
In either case $L_1(s, \pi; r \otimes \eta)$ is defined
to be the gcd of a family of global integrals 
$$
\langle E(f_s), \varphi\rangle_H := \int_{C(\A_F)H(F) \backslash H(\A_F)} 
E(h, f_s) \varphi(h)\omega^{-1}({\rm det}(h)) dh.
\leqno(6.6)
$$
Here $H$ is the reductive group GL$(2)/F$, resp. $R_{K/F}({\rm GL}(2)/K)
\times {\rm GL}(2)/F$, when $m=1$, resp. $m=2$, with center $C$, and
$E(f_s)$ an Eisenstein series on $H(\A_F)$, resp. GSp$(6, \A_F)$ 
associated to a {\it good} section in a representation 
induced from the Borel subgroup. We refer to the papers [HLR] 
and [PS-R] for details.

It is known that $L_1(s, \pi; r \otimes \eta)$
satisfies (FE) and (MC), and also that if $T$ is a finite set of places
containing the archimedean ones and the places where $\pi$ or 
$\eta$ is ramified,
$$
L^T(s, \pi; r \otimes \eta) \, = \, L_1^T(s, \pi; r \otimes \eta).
\leqno(6.7)
$$
So we will be done (for (FE) and (MC)) if we show the following

\medskip

\noindent{\bf Lemma 6.8} \, \it Let $F, \pi, \eta$ be as in Proposition 6.4.
Then the local factors of $L(s, \pi; r \otimes \eta)$ and
$L_1(s, \pi; r \otimes \eta)$ agree at all the places.
\rm

\medskip

{\it Proof of Lemma}. \, Let $v$ be any place
of $F$ which splits in $K$, say as $w, \theta v$, and
$F_v = K_w = K_{\theta w}$. Then the $v$-factor
of $L(s, \pi; r \otimes \eta)$, resp. $L_1(s, \pi; r \otimes \eta)$,
is simply the triple product factor $L(s, \pi_w \times \pi_{\theta w} 
\times \eta_v)$. Similarly for the $\varepsilon$-factors. 
But it was shown in [Ra1], section 4.4, that 
$$
L(s, \pi_w \times \pi_{\theta w} \times \eta_v) 
= L(s, \pi_w \times \pi_{\theta w} \times \eta_v) \quad
{\rm and} \quad 
\varepsilon(s, \pi_w \times \pi_{\theta w} \times \eta_v) 
= \varepsilon(s, \pi_w \times \pi_{\theta w} \times \eta_v).
\leqno(6.9)
$$
So we get the assertion at any such $v$.

Note that any archimedean $v$ splits in $K$ due to our hypothesis that 
$F$ is totally imaginary. Moreover, any $v$ in $S$ splits in $K$ by hypothesis, and
$\pi$ is unramified at places not above $S$.

So we need only prove the assertion at any {\it finite} place $v$ such that (i) there
is a unique place $w$ of $K$ above $v$, \, (ii)
$\eta_v$ is ramified, \, and (iii) $\pi_w$ is unramified. Denote by $\lambda(\pi_w)$, resp.
$\lambda(\eta_v)$, the (non-negative) {\it index of non-temperedness} of $\pi_w$,
resp. $\eta_v$, as in [Ik1]; it is $0$ if the representation is tempered and equals $t > 0$
if it is a complementary series representation defined by the characters $\nu|.|^t, 
\nu|.|^{-t}$ with $\nu$ unitary. One knows that this index is always less than $1/4$ ([GeJ]). (By the recent results of Kim and Shahidi one knows even that it is less than $1/6$, but we do not need this.) Put 
$$
\lambda(\pi_w, \eta_v) \, = \, 2\lambda(\pi_w) + \lambda(\eta_v).
\leqno(6.10)
$$
Suppose $\eta_v$ is non-tempered. Then $\eta_v$ is 
the twist by a unitary character of an unramified representation, and
the truth of the assertion follows from [PS-R], since $\pi_w$ is also unramified. 

So we may assume that $\eta_v$ is tempered, so that 
$$
\lambda(\eta_v) \, = \, 0 \quad {\rm and} \quad \lambda(\pi_w,\eta_v) \, < \, 1/2.
\leqno(6.11)
$$
When $\eta$ is a subquotient of the principal series representation, the assertion then follows, because of (6.11), from Lemma 2.2 of [Ik2]. 

It remains to consider when $\eta_v$ is supercuspidal. Put
$$
\gamma(s, \pi_w; r \otimes \eta_v) \, = \, 
\varepsilon(s, \pi_w; r \otimes \eta_v)\frac
{L(1-s, \pi_w^\vee; r \otimes \eta_v^\vee)}{L(s, \pi_w; r \otimes \eta_v)}.
\leqno(6.12)
$$
Similarly define $\gamma_1(s, \pi_w; r \otimes \eta_v)$.
Since $\pi_w$ is in the principal series, by
applying Prop.5.1 of [Ik1], we get
$$
\gamma(s, \pi_w; r \otimes \eta_v) \, = \, \gamma_1(s, \pi_w; r \otimes \eta_v).
\leqno(6.13)
$$ 
Thanks to (6.11), $L_1(s, \pi_w; r \otimes \eta_v)$, resp. $L(s, \pi_w; r \otimes \eta_v)$,
has no pole in common with $L_1(1-s, \pi_w^\vee; r \otimes \eta_v^\vee)$,
resp. $L(1-s, \pi_w^\vee; r \otimes \eta_v^\vee)$. (In fact, since $\sigma(\eta_v)$
is irreducible, it can be shown that
$L(s, \pi_w; r \otimes \eta_v) = 1$.) Since the $\varepsilon$-factors are invertible, we get the desired equality of $L$ and $\varepsilon$-factors.

\qed 

\medskip

{\it (E)}: \,  
Let $m=1$, and suppose $L(s, \pi; r \otimes \eta)$ has a pole for an idele
class character $\mu$ of $F$. Then, up to replacing $\pi$ 
by $\pi \otimes |.|^{s_0}$ for some $s_0$, we may assume the pole 
to be at $s=1$. Let $S$
be the finite set of places containing the archimedean and ramified places for $\pi$.
Since the local factors have no zeros, the incomplete $L$-function 
$L^S(s, \pi; r \otimes \eta)$  
also has a pole at $s=1$. It is known that the pole must be simple.
Moreover, by Asai's integral representation ([HLR]), the residue at $s=1$ of 
this incomplete
$L$-function is a non-zero multiple of the $\eta$-period ${\mathcal P}_\eta(f)$ 
(see (5.2)).
This means $\pi$ is distinguished, which is ruled out by our hypothesis. Done.

So let $m=2$. Suppose $L(s, \pi; r \otimes \eta)$ has a pole. By Lemma 6.8, we may work with $L_1(s, \pi; r \otimes \eta)$. Again, after replacing $\pi$ by $\pi \otimes |.|^{s_0}$ for suitable $s_0$, we may assume that the pole 
is at $s=1$. Let $\Omega$ denote the central character of $\pi \otimes \eta$, viewed as a cuspidal automorphic representation of $H(\A_F)$. (Recall
that $H = (R_{K/F} {\rm GL}(2)/K) \times {\rm GL}(2)/F$.) By a result of Ikeda (cf. Proposition 2.3 of [Ik1]), $L_1(s, \pi; r \otimes \eta)$ is entire if $\Omega^2$ is non-trivial.

In our case, by hypothesis, $\eta$, and hence its central character, is unramified at the places in $S$, while the square of the restriction of $\omega_\pi$ to $F$ is ramified at some place above $S$. This implies that $\Omega^2$ is ramified, hence non-trivial, at $S$.  Hence we get (E).

\medskip

{\it (BV)}: \, 
Let $S$ be the (finite) set of places of $F$ containing the
archimedean ones and those finite ones ramifiying for $\pi$ or $\eta$. Then the
integral representation of [PS-R] implies the following:
$$
L_1(s, \pi; r \otimes \eta) \prod_{v \in S} 
\frac{\Psi (f_{v,s}; W_v)}{L_1(s, \pi_{v}; r \otimes \eta_v)} \, = \, 
\langle E(f_s), \varphi \otimes \varphi'\rangle_H,
\leqno(6.14)
$$
where $\varphi$ (resp. $\varphi'$) is a cusp form in the space of 
$\pi$ (resp. $\eta$),
$E(f_s)$ is the Sigel Eisenstein series on GSP$(6)/F$ (see [Ra1], sec.3.4, and [Ik1])
associated to a good section
$f_s = \otimes_v f_{v,s}$, $\Psi (f_{v,s}; W_v)$ is, for each $v$, a local integral having 
$L_1 (s, \pi_{v}; r \otimes \eta_v)$ as its gcd for a suitable
$f_{v,s}$ and Whittaker function $W)v$, 
$$
H \,: = \, \{(g,g') \in R_{K/F}{\text GL}(2)/K \times 
{\text GL}(2)/F \, \vert \, {\text det}(g) = {\text det}(g')\}
$$
with center $C$, 
and  
$$
\langle E(f_s), \varphi\rangle_H := \int_{C(\A_F)H(F) \backslash H(\A_F)} 
E(h, f_s) \varphi(h)dh) = \prod_v \Psi (f_{v,s}; W_v).
$$
Using Lemma 3.4.5 of [Ra1] and the hypotheses of this Theorem, we get

\medskip

\noindent{\bf Lemma 6.15}. \quad \it For each $v$, the function 
$\frac{\Psi (f_{v,s}; W_v)}{L_1(s, \pi_{v}; r \otimes \eta_v)}$ is entire and of bounded order, 
for a suitable choice of $f_{v,s}$.
\rm

\medskip

On the other hand, by Proposition 3.4.6 of [Ra1], we know that $E(f_s)$ is a function of bounded order. 
(Analogous results have been established in a very general setting in a recent preprint of W. Muller [Mul]). Since
$\varphi \otimes \varphi'$ vanishes rapidly at infinity, we deduce, using (6.12), that 
$L_1(s, \pi; r \otimes \eta)$ is of bounded order in vertical strips of finite width. 
The same holds then for $L(s, \pi; r \otimes \eta)$ by Lemma 6.8. 
Furthermore, since this $L$-function 
has an Euler product, it is bounded for large positive $\Re(s)$, and hence also for large negative
$\Re(s)$ by the functional equation. Applying the Phragman-Lindel\"of theorem, 
we then conclude the boundedness in vertical strips
of $L(s, \pi; r \otimes \eta)$ as asserted.

Now we are done with the proof of Poposition 6.4.

\qed

\vskip 0.2in

\section
{\bf Proof of Theorem D}

\bigskip

We begin with the following

\medskip

\noindent{\bf Lemma 7.1} \, Let $K/F$ be a quadratic extension of number fields with
$F$ totally imaginary, and
let $\pi$ be a cuspidal automorphic representation of GL$(2, \A_K)$ such that the finite places $w$ of $K$ where $\pi_w$ is ramified are all of degree $1$ over $F$. Then there exists an irreducible, 
admissible, generic representation $\Pi = \otimes_v \Pi_v$
of GL$(4, \A_F)$ such that, for any cuspidal automorphic representation $\eta$ of
G$(m,\A_F)$, $m = 1,2$, and for any place $v$ of $F$, we have
$$
L(s, \Pi_v \times \eta_v) \, = \, L(s, \pi_v; r \otimes \eta_v)
$$
and
$$
\varepsilon(s, \Pi_v \times \eta_v) \, = \, \varepsilon(s, \pi_v; r \otimes \eta_v).
$$
\rm

\medskip

{\it Proof}. \, Let $v$ be any place of $F$. Consider first
the case when $v$ has a unique divisor $w$ in $K$. By hypothesis, $\pi_w$ is in the principal series i,e, it is an isobaric sum $\mu_1 \boxplus \mu_2$. This means $\sigma_w(\pi_w) \simeq \mu_1 \oplus \mu_2$, and
$$
\Lambda^2({\text Ind}_{K_w}^{F_v}(\sigma_w(\pi_w)) \, \simeq \, {\text Ind}_{K_w}^{F_v}(\mu_1)
\otimes {\text Ind}_{K_w}^{F_v}(\mu_2) \oplus \mu_{1,0}\delta_v \oplus \mu_{2,0}\delta_v.
\leqno(7.2)
$$
Here we have used the fact that the determinant of Ind$_{K_w}^{F_v}(\mu_j)$ is the product of the
restriction $\mu_{j,0}$ of $\mu_j$ to $F_v^\ast$ times the quadratic character $\delta_v$ 
associated to $K_w/F_v$. Note also that
$$
{\text Ind}_{K_w}^{F_v}(\mu_1) \otimes {\text Ind}_{K_w}^{F_v}(\mu_2) \,
\simeq \, {\text Ind}_{K_w}^{F_v}(\mu_1\mu_2)
\otimes {\text Ind}_{K_w}^{F_v}(\mu_1(\mu_2 \circ \theta)),
$$
and
$$
{\text Ind}_{K_w}^{F_v}({\rm det}(\sigma_w(\pi_w))) \, \simeq \, {\text Ind}_{K_w}^{F_v}(\mu_1\mu_2).
$$
Then by the definition (3.17) of the Asai representation, we get
$$
As(\sigma_v({\rm Ind}_{K_w}^{F_v}(\pi_w))) \, \simeq \, {\text Ind}_{K_w}^{F_v}(\mu_1(\mu_2 \circ \theta))
\oplus \mu_{1,0}\delta_v \oplus \mu_{2,0}\delta_v.
\leqno(7.3)
$$
Consequently, if we set
$$
\Pi_v \, = \, I_{K_w}^{F_v}(\mu_1(\mu_2 \circ \theta))
\boxplus \mu_{1,0}\delta_v \boxplus \mu_{2,0}\delta_v,
\leqno(7.3)
$$
we will have the properties asserted in the Lemma.

So we may assume that $v$ splits in $K$. Let $w, \theta w$ be the places above $v$.
In this case Ind$_{K_w}^{F_v}(\sigma_w(\pi))$ is just $\sigma_w(\pi) \oplus \sigma_{\theta w}(\pi)$,
which implies in turn that 
$$
As(\sigma_v({\rm Ind}_{K_w}^{F_v}(\pi_w))) \, \simeq \, 
\sigma_w(\pi) \otimes \sigma_{\theta w}(\pi).
\leqno(7.4)
$$
So we may set
$$
\Pi_v \, = \, \pi_w \boxtimes \pi_{\theta w},
\leqno(7.5)
$$
where $\boxtimes$ is the one constructed in [Ra1], 
and this satisfies the asserted properties
relative to any $\eta$. 

\qed

\medskip

\noindent{\bf Proposition 7.6} \, \it Let $(K/F, \pi, \Pi)$ be as in Lemma 7.1. Then $\Pi$ is an isobaric automorphic representation of GL$(4, \A_F)$. 
\rm

\medskip

For any cuspidal automorphic
representation $\eta$ of GL$(m, \A_F)$, $m=1,2$, we will say, following
Piatetski-Shapiro, that the pair $(\Pi, \eta)$ is {\it nice} if the following hold:
\begin{enumerate}
\item[(MC)] \, $L(s, \Pi \times \eta)$ admits a meromorphic 
continuation
to the whole $s$-plane;
\item[(FE)] \, There is a functional equation
$$
L(1-s, \Pi^\vee \times \eta^\vee) \, = \, 
\varepsilon(s, \Pi \times \eta)
L(s, \Pi \times \eta);
$$
\item[(E)] \, $L(s, \Pi \times \eta)$ is entire; \quad and
\item[(BV)] \, $L(s, \Pi \times \eta)$ is bounded in 
vertical strips of 
finite width.
\end{enumerate}
\rm

\medskip

Now we need to
appeal to the following crucial

\medskip

\noindent{\bf Theorem 7.7} (Cogdell - Piatetski-Shapiro [CoPS1]) \, \it 
Let $T$ be a fixed finite set of finite places
of $F$. Let $\beta$ be an irreducible unitary, admissible, generic
representation of GL$(4, \A_F)$. Suppose that for any
cuspidal automorphic representation $\eta$ of GL$(m, \A_F)$,
$m = 1,2$, which is unramified at $T$, the pair $(\beta, \eta)$ is nice. Then $\beta$ is quasi-automorphic, i.e., there is an isobaric automorphic 
representation $\beta_1$ of GL$(4, A_F)$ such that 
$\beta_v \simeq \beta_{1,v}$ at almost all $v$.
\rm

\medskip

{\it Proof of Proposition 7.6}: \, We may assume that $\pi$ is not distinguished. Let $\Pi$ be as in Lemma 7.1. Pick any cuspidal automorphic representation $\eta$ of GL$(m, \A_F)$ for $m \in \{1,2\}$ which is unramified at the set $T$ of finite places where $\Pi$ is ramified. In view of Theorem 7.7, we have to show that $(\Pi, \eta)$ is nice. But in view of Proposition 6.4, it suffices to show that the $L$ and $\varepsilon$-factors of $(\Pi, \eta)$ agree at every place with the corresponding factors of $(\pi; r \otimes \eta)$. This is what was proved in Lemma 7.1, and so we get the quasi-automorphy of $\Pi$. Let $\Pi_1$ be an isobaric automorphic representation of GL$(4, \A_F)$ which is almost everywhere equivalent to $\Pi$. Since the central characters of $\Pi$ and $\Pi'$ agree almost everywhere, they must be equal. Now we make the following

\medskip

\noindent{\bf Claim 7.8} \, \it Let $v$ be any place. Then for any irreducible admissible representation $\beta$ of GL$(m, F_v)$, $m = 1,2$, we have
$$
L(s, \Pi_v \times \beta) \, = \, L(s, \Pi_{1,v} \times \beta)
$$
and
$$
\varepsilon(s, \Pi_v \times \beta) \, = \, \varepsilon(s, \Pi_{1,v} \times \beta).
$$
\rm

\medskip

We know this at all the unramified places and also, by [Ra1], at all the places which split in $K$. So it remains only to prove the assertion at the finite set $S$ of non-split finite places where $\Pi_v$ is ramified. Fix a place $u$ in $S$ and a global atomorphic representation $\eta$ of GL$(m, \A_F)$, $m = 1,2$, with $\eta_u = \beta$. This is clearly possible for $m=1$, and hence for $m=2$ when $\beta$ in the principal series; if $m=2$ and $\beta$ is square-integrable, use the trace formula to constuct such an $\eta$. Choose a global character $\chi$ which is $1$ at $u$ and is {\it highly ramified} at every $v$ in $S - \{u\}$ (relative to $\pi_v, \eta_v$). Now look at any $v \in S - \{u\}$. Since $\Pi_v$ is by construction attached to the $4$-dimensional representation $As(\sigma(\pi_w))$, $w$ being the unique place above $v$, the $L$ and $\varepsilon$-factors of the pair $(\Pi_v, \eta_v \otimes \chi_v)$ are, by the local Langlands correspondence ([HaT], [He]), the same as those of $As(\sigma(\pi_w)) \otimes \sigma(\eta_v) \otimes \chi_v$. Since $\chi_v$ is highly ramified, by [DeH], $L(s, \Pi_v \times \eta_v \otimes \chi_v)$ is $1$, and the dependence of $\varepsilon(s, \Pi_v \times \eta_v \otimes \chi_v)$ on $\Pi_v$ is only through the central character $\omega_{\Pi_v}$. And by [JPSS1], the analogous assertions hold for $L(s, \Pi_{1,v} \times \eta_v \otimes \chi_v)$ and
$\varepsilon(s, \Pi_{1,v} \times \eta_v \otimes \chi_v)$. From this we get, by comparing the functional equations of $L(s, \Pi \times \eta \otimes \chi)$ (cf. Proposition 6.4) and $L(s, \Pi_1 \times \eta \otimes \chi)$, the asserted identity at $u$. Since $u$ was arbitrary in $S$, the claim is proved.  

Finally, Claim 7.8 shows, thanks to a result of Jeff Chen (cf. [Ch], [CoPS2]), that $\Pi_v \simeq \Pi_{1,v}$ at every $v$. So $\Pi$ is isomorphic to the isobaric representation $\Pi_1$. Done with the proof of Proposition 7.6.

\qed

\medskip

It remains to prove {\bf Theorem D in the general case}. Let $(K/F, \pi)$ be arbitrary, but with
$\pi$ not distinguished. We will need to recall the following {\it descent criterion}, which is an 
extension of Proposition 4.2 of [B$\ell$R] to the case of GL$(n)$ for arbitrary $n$:

\medskip

\noindent{\bf Proposition 7.9} ([Ra1], sec. 3.6) \, \it Fix  $n, p \in \N$ with $p$ prime, and a countably infinite set $J$.
Let $F$ a number field, $\{F_j \, | \, j \in J \}$ a family of cyclic
extensions of $F$ with $[F_j : F] = p$, and for each $j \in J$, let $\Pi_j$ be a
cuspidal automorphic representation of GL$(n, \A_{F_j})$. Suppose that, for
all $j, r \in J$, the base changes of $\Pi_j, \Pi_r$ to the compositum
$F_jF_r$ satisfy
$$
(\Pi_j)_{F_jF_r} \, \simeq \, (\Pi_r)_{F_jF_r}.
\leqno(DC)
$$
Then there exists a unique cuspidal automorphic representation $\Pi$ of
GL$(n, \A_F)$ such that
$$
(\Pi)_{F_j} \, \simeq \, \Pi_j, 
$$
for all but a finite number of $j$ in $J$.
\rm

\medskip

First we show how we may restrict to totally imaginary base fields. Indeed, enumerating the prime ideals of $F$ as $P_1, P_2, \cdots$, we may find quadratic extensions $F_j/F$, disjoint from $K/F$, such that for each $j \geq 1$, $F_j$ is totally imaginary {\it and} $P_j$ splits in $F_j$. Put $K_j = KF_j$, which will be a quadratic extension of $F_j$, and let $\pi_j$ denote $\pi_{K_j}$, the base change of $\pi$ to $K_j$. Let $J$ denote the complement in $\N$ of the finite set, possibly empty, of indices $j$ for which either $\pi_j$ is not cuspidal or $\pi$ is
ramified at $P_j$. Associate irreducible, admissible, generic representations $\Pi_j$ to $(\pi_j, K_j/F_j)$ as in Lemma 7.1. Then by Proposition 7.6, each $\Pi_j$ is an isobaric automorphic representation. It is easy to see that the cuspidality condition (of Theorem D) will be satisfied by $\pi_j, r)$ for all but possibly a finite number of $j$. Shrink $J$ by excluding the indices for which this falis. Then by Proposition 4.1, $\Pi_j$ will be cuspidal for each $j$ in $J$. Applying Proposition 7.9, we then get the existence of a common descent $\Pi$ on GL$(4)/F$. By construction, for any $j$ in $J$, the prime $P_j$ splits, say into $Q_j, \overline Q_j$ in $F_j$, and consequently,
$$
L(s, \Pi_{P_j}) = L(s, (\Pi_j)_{Q_j}) = L(s, (\pi_j)_{Q_j};r) = L(s, \pi_{P_j}; r).
\leqno(7.10)
$$
Since this holds at almost all primes $P_j$, $\Pi$ is the desired weak lifting. 

\medskip

Next we reduce to the case when the finite places where $\pi$ is ramified are unramified (for $K/F$) of degree $1$ over $F$. Indeed, let $S$ be the (finite) set of finite places $w$ of $K$ where $\pi$ is ramified and $w$ is either ramified or is of degree $2$ over $F$, so that every such $w$ sits above a unique place $u(w)$ of $F$. Put $S_0 = \{u(w) \vert w \in S\}$, and write $P(w)$ the prime ideal of $F$ defined by $u(w)$. Let $\{P_j \vert j \geq 1\}$ be the set of primes of $F$, and let $J$ denote the complement of the union of $S_0$. For each $j \in J$ choose as above a quadratic extension $F_j$ in which the prime $P_j$ of $F$ splits. We have a ot of freedon in choosing such an $F_j$, and we can take it such that for each $w \in S$, there is a unique place $w(j)$ of $F_j$ such that ${F_j}_{w(j)} \simeq K_w$. Then if we ut $K_j = KF_j$, the base change $\pi_j = \pi_{K_j}$ has the desired ramification property for each $j \in J$. Once we construct $\Pi_j$ over each $F_j$, we can again find a common descent to $F$ using Proposition 7.9. The identity (7.10) will again hold for almost all $j$, and so $\Pi$ will be a weak lifting of $(\pi,r)$.

\medskip

So we may assume from now on that 
\begin{enumerate}
\item[(i)] \, $F$ is totally complex, 
\item[(ii)] \, $\pi$ is a non-distinguished, cuspidal automorphic representation of GL$(2,\A_K)$, and 
\item[(iii)] \, $\pi$ is unramified at the primes of $K$ which are inert or ramified over $F$.
\end{enumerate}

Then the hypotheses of Lemma 7.1 are satisfied, and so we have the existence of the weak lift $\Pi$ on GL$(4)/F$ by Proposition 7.6. It is also the strong lift by Proposition 4.1.
 
Theorem D is now proved.

\qed

\vskip 0.2in

\section{\bf New cases of Artin's conjecture}

\bigskip

Let $\rho, \rho'$ be continuous $\C$-representations of solvable GO$(4)$-type.
By Theorem A, they are modular, associated to isobaric automorphic representations $\pi$, $\pi'$
of GL$(4, \A_F)$. Then 
$$
L(s, \rho \otimes \rho') \, = \, L(s, \pi_f \times \pi'_f).
\leqno(8.1)
$$
The Rankin-Selberg theory of Jacquet, Piatetski-Shapiro and Shalika ([JPSS2], [JS]), 
and of Shahidi ([Sh1,2]), says that the $L$-function 
on the right is entire (see [MW] and the references therein). So Corollary B follows immediately.

\medskip

Now we show how this gives new examples where Artin's conjecture holds. Fix any quadratic extension $E/F$ with non-trivial
automorphism $\theta$. Pick three distinct primes $Q_1, Q_2, Q_3$ of $E$ which are inert over $F$. For each $j \leq 3$, let $P_j$ denote the prime of $F$ below $Q_j$, and $\F(j)$ (resp. $\F_0(j)$) the residue field $\mathfrak O_E/Q_j$ (resp. $\mathfrak O_F/P_j$). Note that $\theta$ induces the non-trivial automorphisms, again denoted by $\theta$, of the $\F(j)$, fixing $\F_0(j)$ pointwise. Choose polynomials $f_j(X) \in \F(j)[X]$, $j \leq 3$, and a polynomial $f(X) \in \mathfrak O_K[X]$, as follows:
\begin{enumerate}
\item[(i)] \, $f_1$ is irreducible of degree $4$, and it does not belong to $\F_0(1)[X]$;
\item[(ii)] \, $f_2$ is the product of an irreducible quadratic polynomial $g_2(X)$ with two distinct linear polynomials;
\item[(iii)] \, $f_3$ is the product of an irreducible cubic polynomial and a linear polynomial;
\item[(iv)] \, $f$ is a quartic polynomial with discriminant $D(f)$, such that $f(X) \equiv f_j(X) ($mod $Q_j)$ for each $j \leq 3$; and
\item[(v)] \, $D(f)^\theta$ is not a square in $E[\sqrt{D(f)}]$. 
\end{enumerate}
(i)-(iv) are easy to achieve. If a particular choice of $f$ satisfying (i)-(iv) fails to satisfy (v), we may then replace it by its sum with a suitable polynomial $h \in \mathfrak O_E[X]$ of degree $\leq 4$ with $h \equiv 0 ($mod $Q_j)$
for each $j$, such that (v) holds. 

\medskip

\noindent{\bf Lemma 8.2} \, \it Let $E/F, f, f_j, j \leq 3$ be as above, satisfying (i)-(v). 
Denote by $K$, resp. $K^\theta$, the splitting field of $f$, resp. $f^\theta$, over $E$. Then $K, K^\theta$ are linearly disjoint over $E$. Moreover, as abstract groups,
$$
{\rm Gal}(K/E) \, \simeq  \, {\rm Gal}(K^\theta/E) \, \simeq \, S_4.
$$
\rm

\medskip

{\it Proof}. \, It is standard that the conditions (i)-(iv) imply that $K/E$ is an $S_4$-extension. This is because these properties realize Gal$(K/E)$ as a transitive subgroup of $S_4$ containing a transposition (coming from $f_2$) and a $3$-cycle (coming from $f_3$), and any such group must be all of $S_4$. Since $f_1$ does not belong to $\F_0(1)[X]$, $f$ does not belong to $\mathfrak O_F[X]$, and so $K$ does not arise by composing $E$ with an $S_4$-extension of $F$. By construction, $f^\theta$ has the same properties, and so $K^\theta/F$ is also an $S_4$-extension.

It remains to show that the field $L: = K \cap K^\theta$, which contains $E$, is just $E$. 

Suppose $L \ne E$. Put $M = KK^\theta$ and $G = $Gal$(M/E)$. Then
$$
G \, \simeq \, {\rm Gal}(K/L) \times {\rm Gal}(K^\theta/L) \, \subset \, S_4 \times S_4.
\leqno(8.3)
$$
Since $K, K^\theta$ are Galois over $E$, so is $L$, and thus
Gal$(K/L)$ and Gal$(K^\theta/L)$ are both proper normal
subgroups of $S_4$. It is well known that the only proper
normal subgroups are $A_4$ and the Kelin group $V$. Either
way, Gal$(K/L)$ must be a subgroup of the alternating
subgroup of Gal$(K/E) \simeq S_4$, and so $L$ contains
$E[\sqrt{D(f)}]$. Similarly, Gal$(K^\theta/L)$ must sit inside
the alternating subgroup of Gal$(K^\theta/E) \simeq S_4$, and
thus $L$ must also contain $E[\sqrt{D(f^\theta)}]$. Then the groups
Gal$(E[\sqrt{D(f)}]/E)$ and Gal$(E[\sqrt{D(f^\theta)}]/E)$ are
both quotients of Gal$(K/E)$. Since $A_4$ is the only
subgroup of $S_4$ of index $2$, we must have $E[\sqrt{D(f)}] =
E[\sqrt{D(f^\theta)}]$. Then $D(f^\theta)$ will be a square in
$E[\sqrt{D(f)}]$, which contradicts (v). Hence the Lemma.

\qed

\medskip

There are clearly an infinite number
of choices for $f$ satisfying such conditions. We can also
vary the inert triple $(Q_1,Q_2,Q_3)$.

\medskip

Since $S_4$ is a subgroup of PGL$(2, \C)$, we get a projective
representation 
$$
\overline \sigma: \, {\rm Gal}(\overline F/E) \,
\rightarrow \, {\rm PGL}(2, \C), 
\leqno(8.4)
$$
such that $\overline F^{{\rm Ker}(\overline \sigma)} = K$.
Its $\theta$-conjugate $\sigma^\theta$ corresponds to
$K^\theta$. 

Fix a lifting 
$$
\sigma: \, {\rm Gal}(\overline F/E) \, \rightarrow
\, {\rm GL}(2, \C)
$$ 
of $\overline \sigma$, which is possible by Tate's theorem
(see section 2). Denote by
$L$ the extension of $F$ corresponding to
ker$(\sigma)$. Clearly, $L \supset K$, and we have a central
extension
$$
1 \, \rightarrow {\rm Gal}(L/K) \, \rightarrow \,
{\rm Gal}(L/E) \, \rightarrow \, {\rm Gal}(K/E) \simeq S_4
\, \rightarrow \, 1.
\leqno(8.5)
$$   

Let us now make precise the structure of Gal$(L/F)$. There
are two double covers, up to isomorphism, of $S_4$, denoted 
$\tilde S_4$ and $\hat S_4$. Gal$(L/K)$ will be a cyclic
subgroup $C$ of order $2m$, for some $m \geq 1$, lying in
the center $Z \simeq \C^\ast$ of GL$(2, \C)$, such that
$$
{\rm Gal}(L/E) \, \simeq \, C \times_{\{\pm I\}} S_4^\ast,
\leqno(8.6)
$$
with
$$
S_4^\ast \, = \, \tilde S_4 \quad {\rm or} \quad \hat S_4.
$$
Here $\times_{\{\pm I\}}$ denotes the direct product with the
common subgroup $\{\pm I\}$ amalgamated. One can identify $\tilde S_4$ with 
GL$(3, \F_3)$, equipped with a natural representation $\pi$ into GL$(2, \C)$. The transpositions 
$\tau$ in $S_4$ lift to elements $\tilde \tau$ of order $2$ in $\tilde S_4$. 
Denote by $\tilde A_4$ the inverse image of $A_4$ in $\tilde S_4$. It is the unique double cover of $A_4$, denoted $2A_4$ in the {\it Atlas}, identifiable with SL$(2, \F_3)$. The (matrix representation of the) group $\hat S_4$ 
can be constructed as follows: Multiply all the elements of $\pi(\tilde S_4)$ outside (resp. inside) $\pi(\tilde A_4)$ by $\sqrt{-1}I$ (resp. $I$). Under this (set-theoretic) map from $\tilde S_4$ to $\hat S_4$, each $\tilde \tau$ goes to an element $\hat \tau$ of order $4$, and this is what distinguishes $\tilde S_4$ from $\hat S_4$. For further reference see [Wa],
pp.9--10.

We also have the $\theta$-conjugate representation
$\sigma^\theta$ of Gal$(\overline F/E)$ into GL$(2, \C)$. The
corresponding field 
$L^\theta$ will contain $K^\theta$ and be linearly disjoint
from $L$ over $E$. Clearly, the Galois group $G$ of
$LL^\theta$ over
$F$ is a non-trivial extension of $\Z/2$ by Gal$(L/E) \times
$Gal$(L^\theta/E)$, and the resulting representation 
$$
\rho: \, {\rm Gal}(\overline F/F) \, \rightarrow \, {\rm
GL}(4, \C)
\leqno(8.7)
$$ of
is evidently irreducible and of GO$(4)$-type. 

Now choose a disjoint quadratic extension $E'/F$ and a quartic polynomial 
$g$ satisfying analogous properties over $E'$. We can arrange, 
in an infinite number of ways, for the resulting 
extension of $F$ to be linearly disjoint from $LL^\theta$. 
Denoting by $\rho'$ the corresponding representation of
Gal$(\overline F/F)$, we see that $\rho \otimes \rho'$ is 
{\it irreducible} and satisfies the Artin conjecture by
Corollary B.

It remains to check that this is not covered by known cases in
lower dimensions. For this it suffices to prove the following

\medskip

\noindent{\bf Proposition 8.8} \, \it Let $\rho, \rho'$ be as
above. Then $\rho, \rho'$ and $\rho \otimes \rho'$ are
primitive representations of Gal$(\overline F/F)$.
\rm

\medskip

By a {\it primitive} representation of Gal$(\overline F/F)$ we
mean a representation which is not induced by a representation
of a proper subgroup. 

\medskip

{\it Proof}.

To begin, it is well known that the irreducible $\sigma$ is
primitive; so is $\sigma^\theta$. 

Suppose $\rho$ is not
primitive. Write 
$$
\rho \, = \, {\rm Ind}_N^F(\tau),
$$
for some representation $\tau$ of Gal$(K/N)$, where
$N \ne F$ is an intermediate field. Note that
$N$ could not be contained in $E$, for otherwise the
restriction of $\rho$ would be reducible, contradicting the
fact that
$$
{\rm Res}_E^F(\rho) \, \simeq \, \sigma \otimes \sigma^\theta,
\leqno(8.9)
$$
which is irreducible as $K, K^\theta$ are linearly disjoint
over $E$. Since $[E:F] = 2$, this implies that Gal$(K/E)$ and
Gal$(K/N)$ generate Gal$(K/F)$. Consequently, by Mackey,
$$
{\rm Res}_E^F(\rho) \, \simeq \, {\rm Ind}_R^E(\tau_0),
\leqno(8.10)
$$
whre $R = E \cap N$ and $\tau_0$ is the restriction of $\tau$
to Gal$(K/R)$. In view of (8.9) and (8.10), it suffices to show
that $\sigma \otimes \sigma^\theta$ is primitive. 
Now we appeal to the following

\medskip

\noindent{\bf Theorem} (Aschbacher [A]) \, \it Let $G_1, G_2$ be 
finite groups and $\pi_i: G_i \rightarrow {\rm GL}(V_i)$ be 
finite-dimensional $\C$-representations which are primitive. 
Then $\pi_1 \otimes \pi_2$ is a primitive representation of $G$.
\rm

\medskip

It may be useful to note that when $G_1, G_2$ are solvable, the proof in [A] 
does not appeal to the classficiation of finite grous.

\medskip

We have already shown that our $\sigma,
\sigma'$ are primitive. Moreover, 
Gal$(KK^\theta/E)$ is isomorpic to Gal$(K/E) \times $Gal$(K^\theta/E)$.
So we may conclude that $\sigma_1 \otimes \sigma_2$ remains primitive. 
In fact this can be verified by a direct,
though laborious,  computation (in our special case), 
which is what we did
originally leading us to pose the general question to
Aschbacher, but now that this primitivity question has been
solved in general, we can do no better than to refer to [A].

Similarly, $\rho'$ is primitive. Applying [A], Theorem 1,
again, we see that $\rho \otimes \rho'$ is also primitive.
Done.

\qed

\vskip 0.2in

\section{\bf The strong Dedekind conjecture: An example with
application}

\bigskip

The {\it Dedekind conjecture} asserts that for any finite
extension $N/F$ of number fields, the zeta function of $F$
divides the zeta function of $N$, i.e., $\zeta_N(s) =
\zeta_F(s)L(s)$ with $L(s)$ entire. If
$\tilde N$ denotes the Galois closure of $N$ over $F$, then
one has
$$
\zeta_N(s) \, = \, \zeta_F(s) L(s, a_{N/F}),
\leqno(9.1)
$$
where $a_{N/F}$ is defined by
$$
a_{N/F} \oplus 1_F \, = \, {\rm Ind}_N^F(1_N).
\leqno(9.2)
$$
Here $1_K$ denotes, for any $K$, the trivial representation
of Gal$(\overline K/K)$. Since $1_F$ does not occur in
$a_{N/F}$, the Artin conjecture for $a_{N/F}$ implies the
Dedekind conjecture for $N/F$. The Dedekind conjecture is known
to be true when $N/F$ is Galois, which follows from the work of Aramata and Brauer, 
and when $\tilde N/F$ is
solvable, proved independently by Uchida and Van der Waall. We
refer to [Mu-R] for a detailed discussion of known results, references and
variants.  We should also note that S. Rallis has a very interesting program for 
studying the ratios $\zeta_N(s)/\zeta_F(s)$ via the 
adjoint $L$-functions of cusp forms $\pi$ on GL$([N:F])$. To elaborate a bit, 
the version of the trace formula due to Jacquet and Zagier ([JZ]) 
suggests that the divisibility of $\zeta_K(s)$ 
by $\zeta_F(s)$ for all commutative algebras $K/F$ of dimension $n$ is equivalent to the divisibility
of $L(s, \pi \times \pi^\vee)$ by $\zeta_F(s)$ for all unitary cuspidal representations $\pi$ of 
GL$(n, \A_F)$ of trivial central character. The adjoint $L$-function of $\pi$ is just the 
ratio $L(s, \pi \times \pi^\vee)/\zeta_F(s)$; so the divisibility in the known case $n=2$ 
can be rederived by using the properties
of the symmetric square $L$-functions of GL$(2)/F$ ([GeJ]). Rallis's idea is to study these via 
certain Eisenstein series on larger groups $G$, and the case $n=3$ has been carried out 
in his intriguing joint paper [JiR], written with
Dihua Jiang, where $G = G_2$.

\medskip

We will say that $N/F$ satisfies the {\it strong Dedekind
Conjecture} if
$a_{N/F}$ is modular, i.e., if there exists an isobaric
automorphic representation of GL$(m,
\A_F)$, $m = [N:F]-1$, with the same $L$-function as that of
$a_{N/F}$. This is known in the following cases:
\begin{enumerate}
\item[(i)] \, $N/F$ is {\it Galois and solvable}; \, and
\item[(ii)] \, $[N:F] \, \leq \, 4$.
\end{enumerate}
The first case is by the
work of Arthur and Clozel ([AC]). In the second case one knows
something stronger, not known in case (i) (unless $K/F$ is
nilpotent), namely that {\it every irreducible occurring in
$a_{N/F}$ is modular}. 

When $N$ is a {\it non-normal cubic} extension of $F$,
$a_{N/F}$ is the unique irreducible $2$-dimensional
representation of Gal$(\tilde N/F) \simeq S_3$; it is
dihedral, induced by either of the non-trivial characters of
the normal subgroup $A_3 = \Z/3$. When
$N$ is a {\it non-normal quartic} extension of $F$, we may
assume that Gal$(\tilde N/F)$ is
$A_4$ or $S_4$, as otherwise the image is nilpotent. Then
$a_{N/F}$ is an irreducible representation induced by a
character of a
$2$-Sylow subgroup $P$. In the $A_4$-case, $P$ is the Klein
group $V$, which is normal in $A_4$, and so the modularity can
be deduced from [AC] once again. But in the
$S_4$-case, $P$ is not normal, and one must appeal to the work
of Jacquet, Piatetski-Shapiro and Shalika ([JPSS1]) on GL$(3)$;
their method is the converse theorem requiring only twists by
characters.

\medskip

\noindent{\bf Theorem 9.3} \, \it Let $F$ be a number
field, $\rho$ a solvable $\C$-representation of
Gal$(\overline F/F)$ of GO$(4)$-type, and $M$ the fixed
field of the kernel of
$\rho$. Then the strong Dedekind conjecture holds for
any intermediate field $N$ (of $M/F$) such that $[N:F] = 3^a$,
$a \geq 0$. In fact, any irreducible summand of $a_{N/F}$ is
modular.
\rm

\medskip

Before beginning the proof, we will indicate some
consequences. Given any cuspidal automorphic representation
$\pi$ of GL$(n, \A_F)$, and a finite extension $N/F$, one can
formally define the {\it base change} $\pi_N$ as an {\it
admissible} representation of GL$(n, \A_N)$. Here is a way to do it.
By the local Langlands correspondence for GL$(n)$ proved by
Harris-Taylor ([HaT]) and Henniart ([He]), we may associate to $\pi_v$, 
at any place $v$ of $F$, a well defined local base change $\pi_{N_w}$ of 
GL$(n, K_w)$, for any place $w$ of $N$ above $v$. To be precise, we take $\pi_{N_w}$
to be the representation associated to the restriction to
$W_{K_w} \times {\rm SL}(2, \C)$ of the $n$-dimensional representation $\sigma(\pi_v)$
of $W_{F_v} \times {\rm SL}(2, \C)$ defined by $\pi_v$. When $\pi_v$ 
is unramified, $\pi_{N_w}$ will
also be unramified, so that the restricted tensor product of the 
$\pi_{N_w}$, as $w$ runs over all the places of $N$, makes sense as an admissible
representation of GL$(n, \A_N)$. Let $\pi_N$ denote $\otimes_w \pi_{N_w}$. Of course it is
a big open problem to know that $\pi_N$ is
automorphic, which is unknown except when $N/F$ is solvable
and {\it normal} ([AC]), and we do not adddress this difficult
question here - at all. But it turns out one can still say a little
bit about the analytic properties of $L(s, \pi_N)$ for special
$N/F$. More precisely, one has the following

\medskip

\noindent{\bf Corollary 9.4} \, \it Let $N/F$ be as in Theorem
9.3. Then for any unitary, cuspidal automorphic representation
$\pi$ of GL$(n, \A_F)$, $L(s, \pi_N)$ admits a meromorphic
continuation to the whole $s$-plane with the expected
functional equation relating $s$ to $1-s$. Moreover, $L(s,
\pi)$ divides $L(s,
\pi_N)$.
\rm

\medskip

Indeed, by Theorem 9.3, $a_{N/F}$ corresponds to an isobaric
automorphic representation $\eta$ of GL$(m, \A_F)$, $m =
[N:F]-1$. It follows that
$$
L(s, \pi_N) \, = \, L(s, \pi)L(s, \pi \times \eta).
\leqno(9.5)
$$
The assertions of the corollary then follow immediately by
applying the Rankin-Selberg theory of Jacquet,
Piatetski-Shapiro and Shalika, and of Shahidi, together with the 
fact that the local Langlands correspondence for GL$(n)$ preserves
$L$ and $\varepsilon$-factors of pairs.

\medskip

\noindent{\bf Remark 9.6} \, We became interested in this
due to the following problem. Consider any elliptic curve
$E$ over
$\Q$. Then for any number field $K$, the Birch and
Swinnerton-Dyer conjecture says that the rank of $E(K)$ as an
abelian group is the order of zero of
$L(s, E_K)$, where $E_K$ denotes the base change of $E$ to
$K$, assuming $L(s, E_K)$ makes sense at $s=1$. Since $E(\Q)
\subset E(K)$, one expects to have 
$$
{\rm ord}_{s=1} L(s, E_K) \, \geq \, {\rm ord}_{s=1} L(s, E).
\leqno(?)
$$
By the monumental work of Wiles ([W]) and Taylor-Wiles ([TW]),
followed by that of Breuil-Conrad-Diamond-Taylor ([BCDT]), we
know that $E$ is modular, i.e., $L(s, E) = L(s-\frac{1}{2},
\pi_f)$, for a unitary cusp form $\pi = \pi_\infty \otimes
\pi_f$ on GL$(2)/\Q$ of weight
$2$, i.e., with $\pi_\infty$ in the lowest discrete series.
So if
$K$ is a $3$-primary extension of $\Q$ contained in the field
cut out by a continuous, solvable representation $\rho$ of
Gal$(\overline \Q/\Q)$ into GO$(4, \C)$, Corollary 9.4 shows
that the expectation $(?)$ does hold in that case.

\medskip

{\it Proof of Theorem 9.3}. \, We will use induction on the
order of $G:= \rho($Gal$(\overline F/F))$. There is nothing to
prove when $|G| = 1$, so assume that $|G| > 1$ and that the
assertion holds for all solvable representations $\rho'$ of
GO$(4)$-type with image of order smaller than $|G|$. 

Write $|G| = 3^{r+a}n$, with $(3,n) =1$, $r \geq 0$. Since $G$
is solvable, there exist subgroups of order $n$, called {\it
Hall subgroups} (relative to the set $S$ of prime divisors of
$n$). They form a single conjugacy class, which we will denote
by
$\mathcal H$. 

Let $R$ denote the image of Gal$(\overline F/N)$ under
$\rho$.  
We claim that $R$ contains some $H$ in
$\mathcal H$. Indeed, since $R$ is solvable of order
$3^rn$, it contains its own Hall subgroups of order $n$,
which must belong to $\mathcal H$. 

Let $L$ be the fixed field
of $H$ so that $K \supset L \supset N \supset F$. 
By the transitivity of induction, we get
$$
a_{L/F} \, \simeq \, {\rm Ind}_N^F({\rm Ind}_L^N(1_K))
\ominus 1_F \, \simeq \, {\rm Ind}_N^F(a_{L/N}) \oplus
a_{N/F}.
\leqno(9.7)
$$
Hence it suffices to prove the assertion for $L/F$.

Suppose $G$ is contained in SGO$(4, \C)$. Then, as we have
seen in section 1, it is given by a product $G_1 \times
G_2$ with each $G_i$ in GL$(2, \C)$. If $H_i$ denotes the
Hall subgroup (relative to $S$) in $G_i$, then $H_1 \times
H_2$ is necessarily in $\mathcal H$. Since $H$ is conjugate
to $H_1 \times H_2$, the corresponding fields have the same
zeta function, and we may, without loss, assume that 
$$
H \, = \, H_1 \times H_2, \quad {\rm and} \quad L \, = \,
L_1L_2,
\leqno(9.8)
$$
where $L_i$ is, for each $i$, the fixed field of $H_i$. 

\medskip

\noindent{\bf Lemma 9.9} \, \it  Every irreducible
subrepresentation $\tau$ of $a_{L_i/F}$, $i \in \{1,2\}$, is of
dimension $\leq 2$.
\rm

\medskip

{\it Proof of Lemma} \, Since $G_i$ is a solvable subgroup of
GL$(2, \C)$, its image $\overline G_i$, say, in PGL$(2, \C)$
is either abelian or dihedral or $A_4$ or $S_4$. It evidently
suffices to consider the latter two cases. Note that $G_i$ is 
an extension of $A_4$ by a {\it central} subgroup $C$, and
the restriction to $C$ of
any irreducible $\tau$ occurring in $a_{L_i/F}$ will be, by
Schur, of the form dim$(\tau)\omega_\tau$, for a character
$\omega_\tau:C \to \C^\ast$. Let $\overline H_i$ denote the
image of the Hall subgroup $H_i$ in $\overline G_i$, and let
$\overline L_i$ be the fixed field of $CH_i$, so that 
$L_i \supset \overline L_i \supset F$.

Suppose first that
$\overline G_i$ is $A_4$. Then
$\overline H_i$ is necessarily the Klein group $V
\simeq \Z/2 \times \Z/2$, which is normal in $A_4$, and the
Galois group of $\overline L_i/F$ is $A_4/V$. 
Clearly,
$a_{\overline L_i/F} \simeq \delta \oplus \delta^2$, where
$\delta$ is a generator of the character group of
Gal$(\overline L_i/F)$. It follows that any irreducible
occurring in $a_{L_i/F}$ is one dimensional, obtained by
pasting onto $\delta^j$, $j \in \{1,2\}$, a character
$\omega$ of $C$ trivial on $\{\pm I\}$.

Now consider when $\overline G_i$ is $S_4$. Now $\overline
H_i$ is a non-normal subgroup of $S_4$ of index $3$. Clearly,
$S_4$ is generated by $A_4$ and $\overline H_i$, and so by Mackey, 
$$
{\rm Res}_{A_4}^{G_i}({\rm Ind}_{\overline H_i}^{G_i}(1)) \,
\simeq \, {\rm Ind}_V^{A_4}(1) \, \simeq \, 1 \oplus \delta
\oplus \delta^2.
$$
Since the cubic character $\delta$ does not extend to 
$S_4$ (whose abelianization is $\Z/2$), it follows that
$$
a_{\overline L_i/F} \, \simeq \, {\rm Ind}_k^F(\delta),
\leqno(9.10)
$$
which is irreducible. Here $k$ denotes the quadratic extension
of $F$ corresponding to (the inverse image of) $A_4$. So any
irreducible $\tau$ occurring in $a_{L-i/F}$ is obtained by
pasting onto this dihedral representation a character $\omega$
of $C$ agreeing on $\{\pm I\}$. The Lemma is now proved.

\qed

\medskip

{\it Proof of Theorem 9.3} (contd.) \, Thanks to (9.8)
we have
$$
{\rm Ind}_L^F(1_L) \, \simeq \, {\rm Ind}_{L_1}^F(1_{L_1})
\otimes {\rm Ind}_{L_2}^F(1_{L_2}).
$$
Hence any irreducible occurring in $a_{L/F}$ is of the
form $\tau_1 \otimes \tau_2$, with $\tau_i$ an irreducible
occurring in $a_{L_i/F}$. By Langlands, $\tau_1$,
resp. $\tau_2$, is modular, associated to a cuspidal
automorphic representation $\pi_1$, resp. $\pi_2$, of GL$(n_i,
\A_F)$, with $n_i = $dim$(\tau_i) \leq 2$. And by [Ra], there
is a cuspidal automorphic representation $\pi_1 \boxtimes
\pi_2$ of GL$(n_1n_2, \A_F)$ having the same $L$-function as
$\tau_1 \otimes \tau_2$. So we are done when $G \subset
$SGO$(4, \C)$.

So we may assume from now on that $G$ is not contained in
SGO$(4, \C)$. Then $G$ has a subgroup $G'$, say, of index $2$
which is a subgroup of SGO$(4, \C)$. Since $H$ is a Hall
subgroup relative to $n$, $H': = H \cap G'$ will necesssarily
be a subgroup of $H$ of index $2$. Then $G$ is generated by
$H$ and $G'$ and so by Mackey,
$$
{\rm Res}_{F'}^{F}({\rm Ind}_{L}^F(1_L)) \,
\simeq \, {\rm Ind}_{L'}^{F'}(1),
$$ 
where $F'$, resp. $L'$, is the quadratic extension of $F$,
resp. $L$, corresponding to $G'$, resp. $H'$. Denote by
$\theta$ the non-trivial automorphism of $F'/F$. 

We have seen that some conjugate of
$H'$ by an element $x$ of $G' \subset G$, is of the form $H_1
\times H_2$, with $H_i$, $i \in \{1,2\}$, being the Hall
subgroup of some solvable $G_i$ in GL$(2, \C)$. So we may,
after replacing $H$ by its conjugate by $x$, assume that
$H'$ is a product group $H_1 \times H_2$. Consequently, for any
irreducible $\tau$ occurring in $a_{L/F}$, we have
$$
{\rm Res}_{F'}^F(\tau) \, \simeq \, \tau_1 \otimes \tau_2,
\leqno(9.11)
$$
where each $\tau_i$ is, by Lemma 9.9, irreducible of
dimension $\leq 2$. The proof of that Lemma shows even that
$\tau_i$ is dihedral if it has dimension $2$. But it is important
to note that $\tau$ can still be primitive, and this fact
provides the content for Theorem 9.3.

Suppose $\tau$ becomes reducible when restricted to
Gal$(\overline F/F')$. Then we must have
$$
\tau \, \simeq \, {\rm Ind}_{F'}^F(\tau'),
\leqno(9.12)
$$
for an irreducible $\tau'$ occurring in $a_{L'/F'}$. If $\pi'$
is the cuspidal automorphic representation of GL$(n', \A_F)$
associated to $\tau'$, with $n' = $dim$(\tau') \leq 2$, then
$\tau$ is modular, associated to the automorphically induced
representation $I_{F'}^F(\pi')$ of GL$(2n', \A_F)$ constructed
by Arthur and Clozel in [AC]. 

Consequently it suffices to consider when $\tau$ remains
irreducible upon restriction to Gal$(\overline F/F')$, being of
the form
$\tau_1 \otimes \tau_2$ (see (9.10)) above). If
either
$\tau_1$ or $\tau_2$ is one dimensional, then $\tau$ has
dimension $\leq 2$, and since its image is solvable, we are
done by taking the Langlands-Tunnell representation
$\sigma(\tau)$ on GL$(2)/F$. 

So we may, and we will,
assume that dim$(\tau_i) = 2$ for each $i$. Now we
may apply Theorem A$^\prime$ and obtain the desired result.
The interesting case is when $\tau_2$ is 
$\tau_1^\theta$, and one gets
$$
\tau \, \simeq \, As(\tau_1),
\leqno(9.13)
$$
and the corresponding automorphic representation $\Pi$ of
GL$(4, \A_F)$ is $As(\sigma(\pi_1))$, as constructed in
Theorem D.

\qed

\vskip 0.2in

\section{\bf Solvable Galois representations of GO$(2m+1)$-type}

\bigskip

In this section we will prove Proposition C.
 
Suppose we are given a continuous irreducible representation of
Gal$(\overline F/F)$ of GO$(n)$-type, wth $n$ odd.
By applying Lemma 1.2 we may, up to replacing $\rho$ by a
one-dimensional twist, which does not affect the conclusion of
the Proposition, assume that the image of $\rho$ lies in O$(n,
\C)$. Let
$K$ be the number field cut out by the kernel of $\rho$ with
(finite) Galois group $G$ over $F$, so that
$\rho$ can be viewed as a faithful representation of $G$. From
the derived series we may extract, by the solvability of $G$,
an elementary abelian 
$p$-group $A$ which is characteristic in $G$, i.e., stable
under any automorphism of $G$. Applying Clifford's theorem, we
see that
$$
{\rm Res}_A^G(\rho) \, \simeq \, m(\chi_1 \oplus \cdots \oplus
\chi_r),
\leqno(10.1)
$$
for some $m,r > 0$ with $mr =n$, and $1$-dimesnional
representations $\chi_1, \cdots, \chi_r$ of $A$ such that
$\chi_i \ne \chi_j$ if $i \ne j$. Moreover, for every $j$
there exists $g_j \in G$ such that 
$$
\chi_j(a) =
\chi_1(g_jag_j^{-1})
\leqno(10.2)
$$ for 
all $a \in A$; hence each $\chi_j$ has the same order, which
must be $p$ as $\rho$ is injective. 

If $p$ is odd, then no
$\chi_j$ is self-dual, while $\rho$ is itself self-dual,
giving a contradiction as $n$ is odd. So $p = 2$. Let 
$$
\rho_1 = m\chi_1 \quad {\rm and} \quad G_1 =
{\rm Stab}_G(\rho_1).
\leqno(10.3)
$$ 
Then
$$
\rho \, \simeq \, {\rm ind}_{G_1}^G(\rho_1)
\leqno(10.4)
$$ 
by Clifford. We
are done if $m=1$.

So we may assume that $m > 1$. If $r = 1$,
$A \simeq \Z/2$, by the faithfulness of $\rho$,
and $\rho(A) = \pm I$. But by construction $A$ is
contained in the commutator subgroup $(G,G)$, which forces
det$(\rho)$ to be trivial on
$A$. On the other hand, since $n$ is odd, det$\rho(A) = -1$, 
resulting in a
contradiction. 

Hence $r$ must be $> 1$ when $m > 1$. Then $G_1$ is a proper
subgroup of $G$ and $\rho_1$ is self-dual by virtue of 
$\chi_1$ being quadratic. Thus
$(\rho_1, G_1)$ satisfies the same hypotheses as $(\rho, G)$.
Since induction is natural in stages, we are done by
infinite descent.

\qed

\vskip 0.2in

\section*{\bf Bibliography}

\begin{description}

\item[{[AC]}] J. Arthur and L. Clozel, {\it Simple Algebras, Base Change
and the Advanced Theory of the Trace Formula}, Ann. Math. Studies {\bf 120}
(1989), Princeton, NJ.

\item[{[A]}] M. Aschbacher, On primitive linear representations of finite groups,
Journal of Algebra {\bf 234}, 627-640 (2000).

\item[{[BaR]}] L. Barthel and D. Ramakrishnan, A non-vanishing result for twists of 
$L$-functions of GL$(n)$, Duke Math. Journal {\bf 74}, no.3, 681-700 (1994).

\item[{[B$\ell$R]}] D. Blasius and D. Ramakrishnan, Maass forms and Galois representations,
in {\it Galois groups over $\Q$}, ed. by Y. Ihara, K. Ribet and J.-P. Serre, MSRI
Publications {\bf 16} (1989). 

\item[{[BCDT]}] C. Breuil, B. Conrad, F. Diamond and R. Taylor, On the modularity of
elliptic curves over $\Q$: wild $3$-adic exercises, to appear in the Journal of the AMS. 

\item[{[Bu]}] J.Buhler, {\it Icosahedral Galois representations}, 
Lecture Notes in Math. {\bf 654}, 
Springer-Verlag (1978).

\item[{[BF]}] D. Bump and S. Friedberg, The exterior square automorphic $L$-functions
on GL$(n)$, Piatetski-Shapiro Festschrift II, Israel Math. Conference Proceedings {\bf 3}, 
47--65 (1990). 

\item[{[BDST]}] K. Buzzard, M. Dickinson, N.I. Shepherd-Barron and R.L. Taylor,
On icosahedral Artin representations, preprint (1999), to appear in the Duke Math Journal

\item[{[BS]}] K.Buzzard and W.Stein, A mod five approach to modularity of icosahedral Galois representations, Preprint (2000), to appear in the Pacific Journal of Math.

\item[{[Ch]}] J.-P. Jeff Chen, {\it Local
factors,  Central characters, and
Representations of the  General Linear Group
over Non-archimedean fields}, Doctoral
dissertation,  Yale University (1996).

\item[{[CoKPSSh]}] J. Cogdell, H. Kim, I. Piatetski-Shapiro and F. Shahidi, 
On lifting from classical groups to GL$_N$, preprint (2000).

\item[{[CoPS1]}] J. Cogdell and I.
Piatetski-Shapiro,  A converse theorem for
GL$(4)$, Math. Research  Letters {\bf 3}
(1996), no. 1, 67-76.

\item[{[CoPS2]}] J. Cogdell and I. Piatetski-Shapiro, Converse theorems for GL$_n$, II, Reine
Angew. Math. {\bf 507}, 165--188 (1999).

\item[{[C]}] L.Clozel, Base change for ${\rm GL}(n)$, Proceedings of the 
International Congress of Mathematicians, 
Vol. 1, 791--797, Amer. Math. Soc., Providence, R.I. (1987)

\item[{[Da]}] E.Dade, Accessible characters are monomial, J. Algebra 
{\bf 117}, no.1, 256--266 (1988).

\item[{[DeH]}] P. Deligne and G. Henniart, Sur la variation, par torsion, 
des constantes locales d'\'equations fonctionnelles de fonctions $L$, Invent. Math. {\bf 64}, 
no. 1, 89--118 (1981).

\item[{[D]}] J. Dieudonn\'e, {\it La g\'eom\'etrie des groupes classiques}, troisi\`eme
edition, Ergebnisse Math., Band {\bf 5}, Springer-Verlag (1971).

\item[{[Fr]}] G. Frey, {\it On Artin's conjecture for odd 2-dimensional representations}, 
Lecture Notes in Math. {\bf 1585}, Springer-Verlag (1994).

\item[{[GeJ]}] S. Gelbart and H. Jacquet, A relation between automorphic
representations of GL$(2)$ and GL$(3)$, Ann. Scient. \'Ec. Norm. Sup. (4)
{\bf 11} (1979), 471--542.

\item[{[GeSh]}] S. Gelbart and F. Shahidi, Boundedness of automorphic 
$L$-functions in vertical strips, preprint (2000).

\item[{[HLR]}] G. Harder, R.P. Langlands and 
M. Rapoport, Algebraische zykeln auf 
Hilbert-Blumenthal-Fl\"achen, Crelles Journal {\bf
366} (1986), 53-120.

\item[{[HaT]}] M. Harris and R. Taylor, On the geometry and cohomology of some 
simple Shimura varieties, preprint (1999), to appear in the Annals of Math. Studies, Princeton.

\item[{[He]}] G. Henniart, Une preuve simple des conjectures de Langlands pour 
${\rm GL}(n)$ sur un corps $p$-adique, Invent. Math. {\bf 139}, no. 2, 439--455
(2000).

\item[{[Ik1]}] T. Ikeda, On the location of poles of the triple
$L$-functions, Compositio Mathematica {\bf 83} (1992), 187-237.

\item[{[Ik2]}] T. Ikeda, On the Gamma factor of the triple $L$-function I, 
Duke Math Journal {\bf 97}, no.2 (99), 301-318.

\item[{[JPSS1]}] H. Jacquet, I. Piatetski-Shapiro and J. Shalika,
Automorphic forms on ${\rm GL}(3)$. II. Ann. of Math. (2) {\bf 109}, 
no. 2, 213--258 (1979). 

\item[{[JPSS2]}] H. Jacquet, I. Piatetski-Shapiro and J. Shalika, 
Rankin-Selberg convolutions, Amer. J. Math.{\bf 105}, 367-464 (983).

\item[{[JS]}] H. Jacquet and J.A. Shalika, Euler products and
the classification of automorphic forms I \& II, Amer. J of
Math. {\bf 103} (1981), 499--558 \& 777--815.

\item[{[JZ]}] H. Jacquet and D. Zagier, Eisenstein series and the Selberg trace formula II, 
Transactions of the AMS {\bf 300} $(1)$, 1-48 (1987).

\item[{[JiR]}] D. Jiang and S. Rallis, Fourier coefficients of Eisenstein series of the
exceptional group of type $G_2$, Pacific Journal of Mathematics {\bf 181}, no.2, 
281-314 (1997).

\item[{[K]}] H. Kim, Functoriality for the exterior square of GL$_4$ and the 
symmetric fourth power of GL$_2$, preprint (2000).

\item[{[KSh]}] H. Kim and F. Shahidi, Functorial products for GL$_2 \times $GL$_3$ 
and functorial symmetric cube for GL$_2$, preprint (2000).

\item[{[La1]}] R.P. Langlands, {\it Base change for GL$(2)$}, Annals of Math. Studies 
{\bf 96}, Princeton (1980).

\item[{[La2]}] R.P. Langlands, On the notion of an automorphic representation. 
A supplement, in {\it Automorphic forms, Representations and $L$-functions}, 
ed. by A. Borel and W. Casselman, Proc. symp. Pure Math {\bf 33}, part 1, 
203-207, AMS. Providence (1979). 

\item[{[La3]}] R.P. Langlands, Letter to A. Weil (January, 1967), posted on
the site
http://www.sunsite.ubc.ca/DigitalMathArchive/Langlands/functoriality.html

\item[{[MW]}] C. Moeglin and J.-L. Waldspurger, Poles des fonctions $L$ de paires pour GL$(N)$, Appendice, Ann. Sci. \'Ecole Norm. Sup. (4) {\bf 22}, 667-674 (1989).

\item[{[Mul]}] W. M\"uller, On the singularities of residual intertwining operators, preprint (2000).

\item[{[MuR]}] M.R. Murty and A Raghuram, Some variations on the Dedekind conjecture, Journal of the Ramanujan Mathematical Society {\bf 15}, No. 4, 225-245 (2000).

\item[{[PS-R]}] I. Piatetski-Shapiro and S. Rallis,
Rankin triple $L$-functions, Compositio Mathematica {\bf 64} (1987), 31-115.

\item[{[Ra1]}] D. Ramakrishnan, Modularity of the Rankin-Selberg $L$-series, and multiplicity one for SL$(2)$, 
Annals of Mathematics {\bf 152}, 43--108 (2000).

\item[{[Ra2]}] D. Ramakrishnan, Arithmetic of Hilbert-Blumenthal surfaces, CMS Conference Proceedings {\bf 7},
285-370 (1987).

\item[{[Se]}] J.-P. Serre, Modular forms of weight one and Galois
representations, in {\it Algebraic Number Fields}, ed. by A. Fr\"ohlich,
Acad. Press (1977), 193-268.

\item[{[Sh1]}] F. Shahidi, On the Ramanujan conjecture and the
finiteness of poles for certain $L$-functions, Ann. of Math. (2) {\bf 127}
(1988), 547--584.

\item[{[Sh2]}] F. Shahidi, A proof of the Langlands conjecture on
Plancherel measures; Complementary series for $p$-adic groups, Ann. of
Math. {\bf 132} (1990), 273-330.

\item[{[T]}] R. Taylor, On icosahedral Artin representations. II, preprint (2000).

\item[{[TW]}] R. Taylor and A. Wiles, Ring-theoretic properties of certain Hecke algebras, 
Ann. of Math. (2) {\bf 141}, no. 3, 553--572 (1995).

\item[{[Tu]}] J. Tunnell, Artin's conjecture for repreesentations of
octahedral type, Bulletin AMS {\bf 5}, no. 2 (1981), 173-175.

\item[{[Wa]}] D. Wales, Quasiprimitive linear groups with quadratic elements, to 
appear in the Journal of Algebra.

\item[{[W]}] A. Wiles, Modular elliptic curves and Fermat's last theorem, Ann. of Math. (2) {\bf 141}, 
no. 3, 443--551 (1995).

\bigskip

\end{description}

\vskip 0.3in

\noindent
Dinakar Ramakrishnan

\noindent
Department of Mathematics

\noindent
California Institute of Technology, Pasadena, CA 91125.

\noindent
dinakar@its.caltech.edu

\end{document}